\input amstex\documentstyle{amsppt}  
\pagewidth{12.5cm}\pageheight{19cm}\magnification\magstep1
\topmatter
\title $\bold Z/m$-graded Lie algebras and perverse sheaves, IV\endtitle
\author George Lusztig and Zhiwei Yun\endauthor
\address{Department of Mathematics, M.I.T., Cambridge, MA 02139}\endaddress
\thanks{The first author was supported in part by NSF grant DMS-1566618.
The second author was supported in part by the Packard Foundation.}\endthanks
\endtopmatter   
\document

\define\tco{\ti{\co}}

\define\dw{\dot w}

     \define\bco{\bar{\co}}

\define\mpb{\medpagebreak}

\define\bN{\bar N}

\define\hW{\hat W}

\define\frl{\forall}

\define\si{\sim}
\define\wt{\widetilde}
\define\sqc{\sqcup}

\define\ovsc{\overset\cir\to}
\define\qua{\quad}

\define\lb{\linebreak}

\define\op{\oplus}
   
\define\part{\partial}
\define\emp{\emptyset}
\define\imp{\implies}
\define\ra{\rangle}
\define\n{\notin}
\define\iy{\infty}
\define\m{\mapsto}
\define\do{\dots}
\define\la{\langle}
\define\bsl{\backslash}

\define\lra{\leftrightarrow}
\define\Lra{\Leftrightarrow}

\define\sm{\smallmatrix}
\define\esm{\endsmallmatrix}
\define\sub{\subset}    
\define\bxt{\boxtimes}
\define\T{\times}
\define\ti{\tilde}
\define\nl{\newline}
\redefine\i{^{-1}}
\define\fra{\frac}
\define\un{\underline}
\define\ov{\overline}
\define\ot{\otimes}
\define\bbq{\bar{\QQ}_l}

\define\ad{\text{\rm ad}}
\define\Ad{\text{\rm Ad}}
\define\Hom{\text{\rm Hom}}

\define\ind{\text{\rm ind}}
\define\Res{\text{\rm Res}}

\define\sg{\text{\rm sgn}}
\define\tr{\text{\rm tr}}

\define\a{\alpha}
\redefine\b{\beta}
\redefine\c{\chi}
\define\g{\gamma}
\redefine\d{\delta}
\define\e{\epsilon}
\define\et{\eta}

\redefine\o{\omega}
\define\p{\pi}
\define\ph{\phi}

\define\r{\rho}
\define\s{\sigma}
\redefine\t{\tau}
\define\th{\theta}

\redefine\l{\lambda}
\define\z{\zeta}
\define\x{\xi}

\define\Om{\Omega}

\define\Th{\Theta}
\redefine\L{\Lambda}

\define\bb{\bold b}

\define\kk{\bold k}

\define\qq{\bold q}

\redefine\ss{\bold s}
\redefine\tt{\bold t}

\define\BB{\bold B}

\define\FF{\bold F}

\define\LL{\bold L}
\define\MM{\bold M}
\define\NN{\bold N}

\define\QQ{\bold Q}
\define\RR{\bold R}

\define\VV{\bold V}

\define\ZZ{\bold Z}
\define\XX{\bold X}
\define\YY{\bold Y}

\define\ca{\Cal A}

\define\cd{\Cal D}
\define\ce{\Cal E}

\define\cg{\Cal G}
\define\ch{\Cal H}
\define\ci{\Cal I}

\define\ck{\Cal K}
\define\cl{\Cal L}
\define\cm{\Cal M}

\define\co{\Cal O}

\define\car{\Cal R}
\define\cs{\Cal S}

\define\cx{\Cal X}

\define\fb{\frak b}

\define\fd{\frak d}
\define\fe{\frak e}
\define\ff{\frak f}
\define\fg{\frak g}

\define\fl{\frak l}
\define\fm{\frak m}

\define\fp{\frak p}
\define\fq{\frak q}

\define\fs{\frak s}
\define\ft{\frak t}
\define\fu{\frak u}

\define\fA{\frak A}

\define\fL{\frak L}

\define\tb{\ti b}

\define\tz{\ti z}

\define\tK{\ti K}

\define\tR{\ti R}

\define\sha{\sharp}

\define\tce{\ti\ce}
\define\bul{\bullet}

\define\che{\check}

\define\cir{\circ}

\define\chR{\che R}    

\define\ALV{A}
\define\BBD{BBD}
\define\BS{BS}
\define\DY{D}
\define\JAC{J}
\define\KLI{KL1}
\define\KLII{KL2}
\define\MADI{L1}
\define\GRE{L2}
\define\ORA{L3}
\define\ICC{L4}
\define\CSI{L5}
\define\CSV{L6}
\define\CAN{L7}
\define\GRAI{L8}
\define\GRAII{L9}
\define\LYI{LY1}
\define\LYII{LY2}
\define\LYIII{LY3}
\define\MA{Ma}
\define\MO{M}
\define\KO{K}
\define\KR{KR}
\define\SHO{Sh1}
\define\SHOGR{Sh2}
\define\SPR{Sp}
\define\VA{Va}
\define\VB{V}
\head Contents\endhead

Introduction

1. $\ZZ/m$-graded root systems.

2. Weyl group representations.

 3. Cyclically graded Lie algebras.

4. Parabolic restriction.

 5. The set $\ci'_0$.

6. Inner product.

7. Induction.

8. Proofs.

Appendix A. An alternative definition of the PBW basis.

Appendix B (by G. Lusztig) $\ZZ$-graded root systems.

\head Introduction\endhead
\subhead 0.1\endsubhead
Let $G$ be a connected reductive group over an algebraically closed field $\kk$ of 
characteristic $p\ge0$ and let $\fg$ be the Lie algebra of $G$. (In the case where $p>0$ we shall
assume that $p$ is a large prime number so that we can operate with Lie algebras as if we were in 
characteristic $0$.) Let $\fg^{nil}$ be the variety of nilpotent elements of $\fg$. We consider 
the adjoint action of $G$ on $\fg^{nil}$; let $G\bsl\fg^{nil}$ be the set of orbits. 

The classification of $G$-orbits on $\fg^{nil}$ has been completed in the 1959 paper of Kostant 
\cite{\KO}. 
Here is some history of this classification. We can assume that $G$ is adjoint simple.
In the case where $G$ is of type $A_n$, the classification was done by Weierstrass 
(1868) and Jordan (1870). In the case where $G$ is of type $B,C$ or $D$ the 
classification was done by 
J. Williamson (1937). Let $J$ be the set of Lie algebra homomorphisms 
$\fs\fl_2@>>>\fg$. Now $G$ acts naturally on $J$; 
let $G\bsl J$ be the set of $G$-orbits on $J$. In 1942, Morozov \cite{\MO} showed
that the map $J@>>>\fg^{nil}$ given by $\ph\m\ph\left(\sm0&1\\0&0\esm\right)$ is 
surjective, hence it induces a surjective map $\th:G\bsl J@>>>G\bsl\fg^{nil}$. (A gap 
in Morozov's proof was filled by Jacobson \cite{\JAC} in 1951.) In 1944, Malcev 
\cite{\MA} showed that $G\bsl J$ is finite; using this and \cite{\MO},\cite{\JAC}, it follows that
$G\bsl\fg^{nil}$ is finite. In 1952, Dynkin \cite{\DY} gave a classification of the 
$G$-orbits on $J$. Finally, in 1959, Kostant \cite{\KO} showed that
$\th:G\bsl J@>>>G\bsl\fg^{nil}$ is injective (hence bijective). This implies a
classification of $G$-orbits on $\fg^{nil}$ (it is the same as the classification of 
$G$-orbits on $J$ which was known from \cite{\DY}). 

\subhead 0.2\endsubhead
Let $T$ be a maximal torus of $G$; let $\ft$ be the Lie algebra of $T$. 
Throughout this paper we assume that $m\in\ZZ_{>0}\cup\{\iy\}$ is given.
If $m<\iy$ we assume that we are given a $\ZZ/m$-grading $\fg=\op_{j\in\ZZ/m}\fg_j$ 
of $\fg$ (see 3.2) such that $\ft\sub\fg_0$. 
If $m=\iy$ we assume that we are given a $\ZZ$-grading $\fg=\op_{N\in\ZZ}\fg_N$ 
of $\fg$ (see B.2) such that $\ft\sub\fg_0$. 
Let $G_0$ be a closed connected subgroup of $G$ whose 
Lie algebra is $\fg_0$. Let $\fg_1^{nil}=\fg_1\cap\fg^{nil}$, a closed subvariety of $\fg_1$ stable under 
the adjoint $G_0$-action. (If $m=\iy$ we have $\fg_1^{nil}=\fg_1$.) The (adjoint)
$G_0$-action on $\fg_1^{nil}$ has only finitely many orbits.
(For $m=2$ this is a result of Kostant and Ralis \cite{\KR}; this was extended to the case $m<\iy$ by Vinberg
\cite{\VB}.)
Let $G_0\bsl\fg_1^{nil}$ be the set of $G_0$-orbits in $\fg_1^{nil}$. Let
$\ci=\ci(\fg_1)$ be the (finite) set of pairs 
$(\co,\ce)$ where $\co\in G_0\bsl\fg_1^{nil}$ and $\ce$ is an irreducible 
$G_0$-equivariant local system on $\co$ (up to isomorphism). For $(\co,\ce)\in\ci$ we 
denote by $\ce^\sha$ the intersection cohomology complex of the closure $\bco$ of $\co$ with 
coefficients in $\ce$, extended by $0$ on $\fg_1^{nil}-\bco$. 
For $(\co,\ce),(\tco,\tce)$ in $\ci$ we define $P_{\tco,\tce;\co,\ce}\in\NN[v\i]$ by
$$P_{\tco,\tce;\co,\ce}=\sum_{a\in\NN}P_{a;\tco,\tce;\co,\ce}v^{-a}$$
where $P_{a;\tco,\tce;\co,\ce}\in\NN$ is the number of times $\tce$ appears in a decomposition of 
the $a$-th cohomology sheaf of $\ce^\sha$ restricted to $\tco$ as a direct sum of irreducible 
local systems and $v$ is an indeterminate. The study of the polynomials 
$P_{\tco,\tce;\co,\ce}$ is of considerable interest. In the case where $m=1$ they appear in
the representation theory of finite reductive groups as certain character values at unipotent
elements; an algorithm for computing them was given in \cite{\CSV}, generalizing earlier work
of the first author \cite{\GRE}, Shoji \cite{\SHOGR} and Beynon-Spaltenstein \cite{\BS}.
In the case where $m=\iy$ they appear in
multiplicity formulas for standard modules of affine Hecke algebras with possibly unequal
parameters; an algorithm for computing them was given in
\cite{\GRAI}, \cite{\GRAII}. In the case where $m<\iy$ they seem to play a role in the character
formulas for double affine Hecke algebras \cite{\VA}, \cite{\LYIII}; 
an algorithm for computing them (except for an 
indexing issue) was given in \cite{\LYI}, \cite{\LYII}. (In these references $G$ is assumed to be
semisimple, simply connected, but for the purposes of this paper these assumptions are not essential.)

In this paper we focus for simplicity on a certain subset $\ci_0=\ci_0(\fg_1)$ (see 3.3) which we call 
the {\it principal block}. In the case where
$m=1$ so that $\fg$ is ungraded and $G_0=G,\fg_1=\fg$, $\ci_0$ is the set of all
$(\co,\ce)$ which appear in the Springer correspondence \cite{\SPR}, see \S5.
We shall consider the square matrix $\MM$ indexed by $\ci_0\T\ci_0$ whose
$((\tco,\tce),(\co,\ce))$ entry is the polynomial 
$$(-v)^{\dim\co-\dim\tco}P_{\tco,\tce;\co,\ce}\in\ZZ[v].$$     
We have the following result (in the case $m=\iy$ this is contained in \cite{\GRAII}):

\proclaim{Theorem 0.3} One can define in a purely combinatorial way a finite set $\BB$ and a 
matrix $\MM'$ of polynomials in $\ZZ[v]$ indexed by $\BB\T\BB$ so that the following holds.
There is an explicit bijection $h:\BB@>\si>>\ci_0$ under which $\MM'$ becomes $\MM$.
\endproclaim

\subhead 0.4\endsubhead
Here ``purely combinatorial'' means that the definition is purely in terms of the root system $R$
of $G$ with respect to $T$ with its $\ZZ/m$-grading (or $\ZZ$-grading)
induced from that of $G$; the group $G$ itself is not 
used in the definition of $\BB$ and $\MM'$. This is reminiscent of the main result of \cite{\KLII}, where 
the polynomials describing the local intersection cohomology of a Schubert variety are identified
with the polynomials of \cite{\KLI} which are defined purely in terms of the root system
(or more precisely the Weyl group). The analogy goes further: in both cases ``semilinear algebra''
(in the form of a bar operator $f\m\bb(f)$, that is the $\QQ$-algebra involution of $\QQ(v)$  
such that $\bb(v^n)=v^{-n}$ for $n\in\ZZ$) plays a key role. 
In our case the set $\BB$ appears as a canonical basis of a $\QQ(v)$-vector space $\VV$ attached 
to the root system with its grading. Following
an idea from \cite{\GRAII} we see that $\VV$ has also another basis $Z$ (which we call the PBW
basis, in analogy with the theory of canonical bases arising from quantum groups \cite{\CAN}) and which is in
natural bijection with $\BB$. Note that
both $\BB$ and $Z$ are defined purely combinatorially, but the
proof that these are well defined is not purely combinatorial, it relies on the geometry of $G$.
(In this respect our results are less satisfactory than those in \cite{\KLI}, \cite{\KLII}.)
The matrix $\MM'$ appears as the transition matrix between the bases $\BB$ and $Z$.

\mpb

Let $\c':\ci_0@>>>G_0\bsl\fg_1^{nil}$ be the map $(\co,\ce)\m\co$. We have the following result.

\proclaim{Theorem 0.5} One can define in a purely combinatorial way a finite set $\Th$ and a
surjective map $\c:\BB@>>>\Th$ so that the following holds. There is an explicit bijection
$h':\Th@>\si>>G_0\bsl\fg_1^{nil}$ such that $\c'h=h'\c:\BB@>>>G_0\bsl\fg_1^{nil}$.
\endproclaim
In fact $\Th$ appears as a certain finite set of facets  (which we call rigid)
of an affine  hyperplane arrangement associated to 
$R$ (with its grading) modulo the action of a certain subgroup of the Weyl group.
In the case where $m=1$, this hyperplane arrangement is the standard one associated to the affine Weyl group coming from $R$.
If $m=\iy$ the hyperplanes in the arrangement all pass through $0$. 

\mpb

We now state two results about the fibres of the map $\c:\BB@>>>\Th$.

\proclaim{Theorem 0.6} For any $\o\in\Th$ one can define in a purely combinatorial way
a certain set $\hW^{[\o]}$ of irreducible representations of a certain Weyl group (depending on $\o$)
and a canonical bijection $\c\i(\o)\lra\hW^{[\o]}$.
\endproclaim

\proclaim{Theorem 0.7} Assume that $m=1$. Let $\hW$ be the set of isomorphism classes of irreducible representations
(over $\QQ$) of the Weyl group of our root system. One can define in a purely combinatorial way a partition
$\hW=\sqc_{\o\in\Th}\hW^\o$ and, for any $\o\in\Th$, a canonical bijection $\c\i(\o)\lra\hW^\o$.
\endproclaim
This is essentially the same as the Springer correspondence \cite{\SPR} except that, unlike our bijection,
the Springer correspondence is not purely combinatorial; its definition is based on geometry.
Note also that $\hW^{[\o]}$ in 0.6 (with $m=1$) is not necessarily the same as $\hW^{\o}$ in 0.7, although the two are
in canonical bijection, see 8.7.

\subhead 0.8\endsubhead
The proof of each of the Theorems 0.3, 0.5, 0.6, 0.7 relies in part on
the semilinear algebra computations in the $\ZZ$-graded case given in \cite{\GRAII}.
But this goes also in the opposite direction: the proof of Theorem 0.6 for $m=\iy$
requires arguments in the case $m=1$. 

\subhead 0.9\endsubhead
The paper has two appendices. In Appendix A we give a definition of $\BB$ and $Z$ when $m<\iy$ which
does not rely on the results in \cite{\GRAII}; the definition of $\BB$
is a simplification of one in \cite{\LYII}. This gives another proof of Theorem 0.3 which is not
relying on \cite{\GRAII} (it still relies on \cite{\GRAI}). But this approach is not capable
of proving Theorems 0.5, 0.6, 0.7.

In Appendix B we reformulate the results in \cite{\GRAII} in a form that can be used in this paper.

\subhead 0.10. Notation\endsubhead
If $A$ is a subset of a vector space $E$ we denote by $\la A\ra$ the subspace of $V$
generated by $A$.

Let $\ca=\ZZ[v,v\i]$.

If $V,V'$ are $\QQ(v)$-vector spaces, a $\QQ$-linear map $\b:V@>>>V'$ is said to be 
{\it semilinear} if $\b(fx)=\bb(f)\b(x)$ for any $x\in V,f\in\QQ(v)$.

If $x\in\QQ-\{0\}$ we set $\sg(x)=1$ if $x>0$, $\sg(x)=-1$ if $x<0$.

For any linear algebraic group $\cg$ let $\fL\cg$ be the Lie algebra of $\cg$. 

All algebraic varieties are assumed to be over $\kk$. For an algebraic variety $X$ we denote by
$\cd(X)$ the bounded derived category of $\bbq$-sheaves on $X$; here $l$ is a fixed prime number
invertible in $\kk$. For $K\in\cd(X)$ let ${}^pH^jK$ be the $j$-th perverse cohomology sheaf of $K$
and let $\ch^jK$ be the $j$-th cohomology sheaf of $K$.

If $F:X@>>>X$ is a map of sets, we write $X^F=\{x\in X;F(x)=x\}$.

If $m\in\ZZ_{>0}$ we denote by $\bN$ the image of $N\in\ZZ$ in $\ZZ/m$; for 
$j\in\ZZ/m$ let $\wt{j}=\{N\in\ZZ;\bN=j\}\sub\ZZ$. 

\head 1. $\ZZ/m$-graded root systems\endhead
\subhead 1.1\endsubhead
In this section we state our main results purely combinatorially in terms of $\ZZ/m$-graded root
systems (with $m<\iy$). Let

(a) $(\YY,\XX,(,),\chR\lra R)$
\nl
be a root system. (We often write $R$ instead of (a).)
Thus, $\YY,\XX$ are $\QQ$-vector spaces of finite dimension, 
$(,):\YY\T\XX@>>>\QQ$ is a perfect pairing and $\chR\sub\YY,R\sub\XX$. Let $W$ be the 
Weyl group of $R$ viewed as a subgroup of $GL(\YY)$ and as a subgroup of $GL(\XX)$.

We view $\YY$ and its subsets with the topology induced from the standard
topology of $\RR\ot_\QQ\YY$.

\subhead 1.2\endsubhead
In this section, until the end of 1.12, we assume that $m<\iy$.
A $\ZZ/m$-grading for $R$ is a collection $(R_j)_{j\in\ZZ/m}$ where $R_j$ are subsets of 
$R$ such that $R=\sqc_{j\in\ZZ/m}R_j$ and such that for $\a\in R_j,\a'\in R_{j'}$ we have 
$\a+\a'\in R\imp\a+\a'\in R_{j+j'}$ and $\a+\a'=0\imp j+j'=0$. 
We assume that a $\ZZ/m$-grading for $R$ is fixed. 
Let $\chR_0$ be the image of $R_0$ under $\chR\lra R$; then
$(\YY,\XX,(,),\chR_0\lra R_0)$ is a root system. Its Weyl group $W_0$ is the subgroup 
of $W$ generated by the reflections with respect to roots in $R_0$. The obvious $W_0$-action 
on $R$ leaves stable each of the subsets $R_j,j\in\ZZ/m$. 

Let $e_{W_0}=\sum_{w\in W_0}v^{2|w|}$ where $w@>>>|w|$ is the length function on $W_0$ for a
Coxeter group structure on $W_0$ determined by any choice of simple roots for $R_0$.
Let $\cs$ be the collection of affine hyperplanes
$$\{\{y\in\YY;(y,\a)=N\}; N\in\ZZ,\a\in R_{\bN}\}.$$
Let 
$$\YY'=\YY-\cup_{H\in\cs}H=\YY-\cup_{j\in\ZZ/m,\a\in R_j}\{y\in\YY;(y,\a)\in\wt{j}\}.$$   
The facets determined by $\cs$ are called {\it $m$-facets}. They can be described as follows.
For $y_1,y_2$ in $\YY$ we write $y_1\si_my_2$ if for any $j\in\ZZ/m$, any $\a\in R_j$
and any $N\in\wt{j}$ we have $(y_1,\a)\ge N\Lra(y_2,\a)\ge N$. If $y_1\si_my_2$ and  
$j\in\ZZ/m$, $\a\in R_j$, $N\in\wt{j}$, then $(y_1,\a)>N\Lra(y_2,\a)>N$. (Indeed, assume 
that $(y_1,\a)>N,(y_2,\a)\not>N$. We must have $(y_2,\a)=N$. We have $-\a\in R_{-j}$, 
$(y_2,-\a)=-N$ hence $(y_1,-\a)\ge-N$ and $(y_1,\a)\le N$, contradicting $(y_1,\a)>N$.)
We deduce that, if $y_1\si_my_2$ and $j\in\ZZ/m$, $\a\in R_j$, $N\in\wt{j}$, then 
$(y_1,\a)=N\Lra(y_2,\a)=N$. Now $\si_m$ is an equivalence relation; the equivalence 
classes are the $m$-facets.

An $m$-facet is said to be an $m$-{\it alcove} if
it is contained in $\YY'$. Let $\un\YY'$ be the set of $m$-alcoves.

For $y,y'$ in $\YY'$ we define
$$\align&\t(y,y')=\t(y',y)\\&=
\sha\{\a\in R_1;((y,\a)-1)((y',\a)-1)<0\}-\sha\{\a\in R_0;(y,\a)(y',\a)<0\}\in\ZZ,\endalign$$
$$(y:y')=e_{W_0}\sum_{w\in W_0}v^{\t(y,w(y'))}\in\ca.$$
Let $\VV'=\VV'_R$ be the $\QQ(v)$-vector space with basis $\{I_\g;\g\in\un\YY'\}$.
\footnote{In \S3, the elements $I_\g$ are interpreted as spiral inductions.}
We 
define a $\QQ(v)$-bilinear form $(:):\VV'\T\VV'@>>>\QQ(v)$ by $(I_\g:I_{\g'})=(y:y')$
where $\g\in\un\YY',\g'\in\un\YY'$, $y\in\g,y'\in\g'$; this is independent of the 
choice of $y,y'$. This form is symmetric since $\t(y,w(y'))=\t(y',w\i(y))$ for $w\in W_0$, 
$y\in\YY',y'\in\YY'$. Let \lb $\car=\{\x\in\VV';(x:\VV')=0\}$, $\VV=\VV_R=\VV'/\car$; now $(:)$ 
induces a nondegenerate symmetric bilinear form $\VV\T\VV@>>>\QQ(v)$ denoted again by $(:)$.

For $\g\in\un\YY'$, the image in $\VV$ of $I_\g\in\VV'$ is denoted again by $I_\g$.
One can show (see 8.3): 

(a) {\it There is a unique semilinear map $\b:\VV@>>>\VV$ such that $\b(I_\g)=I_\g$ for 
any $\g\in\un\YY'$.}

\subhead 1.3\endsubhead
For any $m$-facet $\r$ and any $N\in\ZZ$ we set $R(\r)_N=\{a\in R_{\bN};(y,\a)=N\}$ 
where $y\in\r$; this is independent of the choice of $y$. We set 
$$R(\r)=\sqc_{N\in\ZZ}R(\r)_N=\cup_{j\in\ZZ/m}\{\a\in R_j;(y,\a)\in\wt{j}\}$$
where $y\in\r$. Let $\chR(\r)$ be the image of $R(\r)$ under $\chR\lra R$. Then\lb
$(\YY,\XX,(,),\chR(\r)\lra R(\r))$ is a root system with a $\ZZ$-grading 
(in the sense of B.2) $R(\r)_*=(R(\r)_N)_{N\in\ZZ}$. As in B.2, there is a unique 
element $y_{R(\r)_*}\in\la\chR(\r)\ra\sub\YY$ such that for any $N\in\ZZ$,
$\a\in R_N(\r)$, we have $(y_{R(\r)_*},\a)=N$. We say that $\r$ is {\it rigid} if 
$y_{R(\r)_*}\in\r$ and if $R(\r)_*$ is rigid in the sense of B.7. Let $\un\YY^\bul$ 
be the set of rigid $m$-facets. Now the obvious $W_0$-action on the set of $m$-facets 
preserves the set $\un\YY^\bul$. Let $\un{\un\YY}^\bul$ be the set of $W_0$-orbits on 
$\un\YY^\bul$. 

\subhead 1.4\endsubhead
We now return to a general $m$-facet $\r$. Let 
$$\YY'(\r)=\{y'\in\YY;(y',\a)\ne0\qua\frl\a\in R(\r)\}.$$
On $\YY'(\r)$ we have an equivalence relation where $y_1,y_2$ are equivalent if \lb
$(y_1,\a)(y_2,\a)>0$ for any $\a\in R(\r)$; let $\un{\YY'(\r)}$ be the set of 
equivalence classes, that is, the set of $\iy$-alcoves of $(\YY,\XX,(,),\chR(\r)\lra R(\r))$ 
(as in B.1). We define a map $f_\r:\un{\YY'(\r)}@>>>\un\YY'$ as follows. Let
$y\in\r$ and let $r\in\ZZ_{>0}$ be such that $(y,\a)\in(1/r)\ZZ$ for any $\a\in R$. Let 
$\g\in\un{\YY'(\r)}$ and let $y_1\in\g$. We can assume that $-1<(y_1,\a)<1$ for any
$\a\in R$. We show:

(a) $y+r\i y_1\in\YY'$.
\nl
Assume that $j\in\ZZ/m$, $\a\in R_j$ and $(y+r\i y_1,\a)\in\wt{j}$. Since 
$(y,\a)\in(1/r)\ZZ$, we have $(r\i y_1,\a)=(y+r\i y_1,\a)-(y,\a)\in(1/r)\ZZ$.
Combining this with $-1/r<(r\i y_1,\a)<1/r$, we see that $(r\i y_1,\a)=0$ hence 
$(y_1,\a)=0$ and $(y,\a)\in\wt{j}$. But $(y_1,\a)=0$ implies $\a\n R(\r)$ hence 
$(y,\a)\n\wt{j}$. This is a contradiction; (a) is proved.

Now let $y'\in\r$ and let $r'\in\ZZ_{>0}$ be such that $(y',\a)\in(1/r')\ZZ$ for any
$\a\in R$; let $y'_1\in\g$ be such that $-1<(y'_1,\a)<1$ for any $\a\in R$. By (a)
we have $y'+(r')\i y_1\in\YY'$. We show:

(b) $y+r\i y_1\si_m y'+(r')\i y'_1$.
\nl
Assume that for some $N\in\ZZ$ and some $\a\in R_{\bN}$, 

$(y+r\i y_1,\a)-N$, $(y'+(r')\i y'_1,\a)-N$ 
\nl
have different signs. If $\a\in R(\r)$, then $(y,\a)=M$, 
$(y',\a)=M'$ for some $M\in N+m\ZZ$, $M'\in N+m\ZZ$; since $y\si_my'$, we have $M=M'$, 
so that 

$M-N+(r\i y_1,\a)$, $M-N+((r')\i y'_1,\a)$ 
\nl
have different signs; since 

$-1/r<(r\i y_1,\a)<1/r$, $-1/r'<((r')\i y'_1,\a)<1/r'$,
\nl
it follows that $M=N$ and that $(r\i y_1,\a)$, $((r')\i y'_1,\a)$ have different signs and 
$(y_1,\a),(y'_1,\a)$ have different signs, contradicting that $y_1\in\g,y'_1\in\g$.

If $\a\n R(\r)$ then $(y,\a)\in(1/r)\ZZ$, $(y,\a)\n N+m\ZZ$ hence $|(y,\a)-N|\ge1/r$; 
since $-1/r<(r\i y_1,\a)<1/r$, we see that $(y,\a)-N$ has the same sign as 
$(y+r\i y_1,\a)-N$. Similarly $(y',\a)-N$ has the same sign as $(y'+(r')\i y'_1,\a)-N$.
Thus, $(y,\a)-N$, $(y',\a)-N$ have different signs. This 
contradicts the fact that $y\si_my'$ and proves (b).

We see that $\g\m y+r\i y_1$ induces a well defined map 
$$f_\r:\un{\YY'(\r)}@>>>\un\YY'.$$

\subhead 1.5\endsubhead
Let $\r$ be an $m$-facet. Let $\VV'_{R(\r)},\VV_{R(\r)},(:)$ be the analogues of $\VV',\VV,(:)$ 
in B.3 when $R_*$ is replaced by the $\ZZ$-graded root system $R(\r)_*$.
 We define a $\QQ(v)$-linear map 
$\VV'_{R(\r)}@>>>\VV'$ by sending the basis element indexed by $\g\in\un{\YY'(\r)}$ to 
$I_{f_\r(\g)}$. One can show (see 8.3):

(a) {\it this maps the radical of $(:)$ on $\VV'_{R(\r)}$ into the radical of $(:)$ on 
$\VV'$ hence it induces a linear map $\VV_{R(\r)}@>>>\VV$ denoted again by $f_\r$.}

\subhead 1.6\endsubhead
We now assume that $\r$ is a rigid $m$-facet. Let ${}^1Z(\r)_{R(\r)}$ be the PBW basis of 
$\VV_{R(\r)}$ defined as in B.6 in terms of $R(\r)_*$ with $\d=1$. Let 
$$[0]_\r=\{y\in\YY;(y,\a)=0\qua\frl\a\in R(\r)\}.$$
Then the subset ${}^1Z(\r)^{[0]_\r}_{R(\r)}$ of ${}^1Z(\r)_{R(\r)}$ defined as in B.6 is nonempty. Let 
$Z^\r=Z^\r_R=f_\r({}^1Z(\r)^{[0]_\r}_{R(\r)})\sub\VV$;
it depends only on the $W_0$-orbit $(\r)$ of $\r$; we
shall write $Z^{(\r)}=Z^{(\r)}_R$ instead of $Z^\r_R$. One can show (see 8.3):

(a) {\it $f_\r$ is a bijection ${}^1Z(\r)^{[0]_\r}_{R(\r)}@>\si>>Z^{(\r)}_R$.}

(b) {\it The union $Z=Z_R:=\cup_{\o\in\un{\un\YY}^\bul}Z^\o_R$ is disjoint.}

(c) {\it $Z_R$ is a basis of the vector space $\VV$ (called a PBW-basis). Let $\cl$ be the
$\ZZ[v]$-submodule of $\VV$ generated by $Z_R$.}

(d) {\it For each $\x\in Z_R$ there is a unique element $\un\x\in\cl$ such that $\un\x-\x\in v\cl$ 
and $\b(\un\x)=\un\x$.}

(e) {\it The map $\x\m\un\x$ is a bijection of $Z_R$ onto a $\ZZ[v]$-basis
$\BB=\BB_R$ of $\cl$ which is
also a $\QQ(v)$-basis of $\VV$ called the canonical basis of $\VV$. Under this bijection,
the subset $Z^\o_R$ of $Z_R$ (where $\o\in\un{\un\YY}^\bul$)
corresponds to a subset $\BB^\o$ of $\BB$; we 
have $\BB=\sqc_{\o\in\un{\un\YY}^\bul}\BB^\o$.}
\nl
Let $\o\in\un{\un\YY}^\bul$. We set
$$d(\o)=\sha\{\a\in R_0;(y',\a)<0\}+\sha\{\a\in R_1;(y',\a)\ge1\}$$
where $\r\in\o$ and $y'\in\r$. One can show (see 8.3):

(f) {\it For any $\x\in Z^\o_R$, $\un\x-\x$ is a linear combination with coefficients in
$v\ZZ[v]$ of elements $\x'\in Z^{\o'}$ where $\o'\in\un{\un\YY}^\bul$ satisfies $d(\o')<d(\o)$.}

(g) {\it Let $\x\in Z^\o,\x'\in Z^{\o'}$ with $\o,\o'$ in $\un{\un\YY}^\bul$. If $\x=\x'$, then
$(\x:\x')\in 1+v\ZZ[v]$. If $\x\ne\x'$, then $(\x:\x')\in v\ZZ[v]$. Moreover, if $\o\ne\o'$ then 
$(\x:\x')=0$.}

(h) {\it Let $\et\in\BB,\et'\in\BB$. If $\et=\et'$, we have $(\et:\et')\in1+v\ZZ[v]$. 
If $\et\ne\et'$ we have $(b:b')\in v\ZZ[v]$.}

\subhead 1.7\endsubhead
Until the end of 1.12 we assume that $m=1$. In this case we have 
$\ZZ/m=\{0\}$ and $R=R_0$, $W=W_0$. We now have 
$$\YY'=\YY-\cup_{\a\in R}\{y\in\YY;(y,\a)\in\ZZ\}.$$
For $y_1,y_2$ in $\YY$ we have $y_1\si_1y_2$ if for any $\a\in R$ and any $N\in\ZZ$ 
we have $(y_1,\a)\ge N\Lra(y_2,\a)\ge N$. For $y,y'$ in $\YY'$ we have
$$\align&\t(y,y')=\t(y',y)\\&
=\sha\{\a\in R;((y,\a)-1)((y',\a)-1)<0\}-\sha\{\a\in R;(y,\a)(y',\a)<0\}\in\ZZ,\endalign$$
$$(y:y')=e_W\sum_{w\in W}v^{\t(y,w(y'))}\in\ca.$$
For any $1$-facet $\r$ and any $N\in\ZZ$ we have $R(\r)_N=\{a\in R;(y,\a)=N\}$ where 
$y\in\r$. We have
$$R(\r)=\sqc_{N\in\ZZ}R_N(\r)=\{\a\in R;(y,\a)\in\ZZ\}$$
where $y\in\r$. In this case, for $\o\in\un{\un\YY}^\bul$ we have
$$d(\o)=\sha\{\a\in R;(y,\a)<0\}+\sha\{\a\in R;(y,\a)\ge1\}$$
where $\r\in\o$ and $y\in\r$. 

\subhead 1.8. Examples\endsubhead
Recall that $m=1$. 

Assume first that $\YY=\XX=\QQ$, $(,)$ is given by the product in $\QQ$, 
$R=\{-1,1\}$, $\chR=\{-2,2\}$. The $1$-alcoves are the subsets $\{x\in\QQ;n<x<n+1\}$ for 
various $n\in\ZZ$. Let $c=\{x\in\QQ;x<-1\}\cup\{x\in\QQ;x>1\}$, $c^0=\{x\in\QQ;-1<x<1\}$. 
For two $1$-alcoves $\g,\g'$ we have in $\VV$:

$(I_\g:I_{\g'})=(1+v^2)(1+v^{-2})$ if $\g\sub c^0$, $\g'\sub c^0$

$(I_\g:I_{\g'})=(1+v^2)(v\i+v)$, if $\g\sub c^0,\g'\sub c$,

$(I_\g:I_{\g'})=2(1+v^2)$ if $\g\sub c$, $\g'\sub c$.
\nl
The canonical basis of $\VV$ is $\{A,B-A\}$ where $A=\fra{1}{v+v\i}I_\g$ with 
$\g\sub c^0$, $B=I_\g$ with $\g\sub c$. We have $(A:A)=1$, $(B-A:B-A)=1$, $(A:B-A)=v^2$. 
The PBW basis is $\{A,-v^2A+(B-A)\}$.

\mpb

Next we assume that $\YY$ has basis $\che\a_1,\che\a_2$,
$\XX$ has basis $\a_1,\a_2$, $(,)$ is given by $(\che\a_i,\a_i)=2$ for $i=1,2$, 
$(\che\a_i,\a_j)=-1$ if $i\ne j$;
$\chR$ consists of $\pm\che\a_1,\pm\che\a_2,\pm(\che\a_1+\che\a_2)$;
$R$ consists of $\pm\a_1,\pm\a_2,\pm(\a_1+\a_2)$.
Let $\g_0,\g_1,\g_2$ be the $1$-alcoves containing $(\che\a_1+\che\a_2)/3$,
$2(\che\a_1+\che\a_2)/3$, $4(\che\a_1+\che\a_2)/3$ respectively.
Then $\{I_{\g_k};k=0,1,2\}$ is a $\QQ(v)$-basis of $\VV$.
We have $(I_{\g_0}:I_{\g_0})=v^{-6}e_W^2$, $(I_{\g_0}:I_{\g_1})=v^{-5}e_W^2$,
$(I_{\g_0}:I_{\g_2})=v^{-3}e_W^2$, $(I_{\g_1}:I_{\g_1})=(v^{-4}+2v^{-2}+3)e_W$,
$(I_{\g_1}:I_{\g_2})=(v^{-2}+4+v^2)e_W$, $(I_{\g_2}:I_{\g_2})=6e_W$.
The canonical basis $\{b_0,b_1,b_2\}$ of $\VV$ satisfies
$I_0=(e_W v^{-3})b_0$, $I_1=b_1+(v\i+v)^2b_0$, $I_2=b_2+2b_1+b_0$.
The PBW basis $z_0,z_1,z_2$ of $\VV$ satisfies
$b_0=z_0$, $b_1=z_1+(v^4+v^2)z_0$, $b_2=z_2+v^2z_1+v^6z_0$.

\subhead 1.9\endsubhead
Recall that $m=1$. Let $\r$ be an $\iy$-facet  (as in B.1).
Let $R'=\{\a\in R;(y',\a)=0\}$ where $y'\in\r$. Let $\chR'$ be the image of $R'$ under 
$\chR\lra R$. Then $(\YY,\XX,(,),\chR'\lra R')$ is a root system.
Let $W'$ be the Weyl group of $R'$, viewed as a subgroup of $W$. Let $\fe$ be a subset of 
$W$ such that $W=W'\fe$, $\sha\fe=\sha W/\sha W'$. For $\e\in\fe$ we set
$$f'_\e=-\sum_{\a\in R';(y',\a)\ne0,(\e(y'),\a)<0\text{ or }(\e(y'),\a)>1}\sg(y',\a).\tag a$$
For any $1$-alcove $\g$ relative to $R$ and any $\e\in\fe$,
$\e(\g)$ is contained in a unique $1$-alcove $\wt{\e(\g)}$ relative to $R'$.
We define $\VV_{R'}$ in terms of $R'$ in the same way as $\VV$ was defined in terms of $R$.

For any $1$-alcove $\g'$ relative to $R'$ we define $I_{\g'}\in\VV_{R'}$ in the same way as 
$I_\g\in\VV$ was defined in terms of $R$. One can show (see 8.3):

(b) {\it There is a unique $\QQ(v)$-linear map $\Res_\r:\VV@>>>\VV_{R'}$ such that for any
$1$-alcove $\g$ relative to $R$ we have }
$$\Res_\r(I_\g)=\sum_{\e\in\fe}v^{f'_\e}I_{\wt{\e(\g)}}.$$

\subhead 1.10\endsubhead
We preserve the setup of 1.9. 
Let $Z_{R'}$ (resp. $\BB_{R'}$) be the PBW basis (resp. canonical basis) of
$\VV_{R'}$ (with $m=1$). Let $\cl_{R'}$ be the $\ZZ[v]$-submodule of $\VV_{R'}$ spanned by 
$Z_{R'}$. One can show (see 8.3): 

(a) {\it Let $\x\in Z$. We have $\Res_\r(\x)=\sum_{\x'\in Z_{R'}}c_{\x',\x}\x'\mod v\cl_{R'}$ 
where $c_{\x',\x}\in\ZZ$.}
\nl
From (a) and 1.6(d),(e) we deduce:

(b) {\it Let $\et\in\BB$. We have $\Res_\r(\et)=\sum_{\et'\in\BB_{R'}}c_{\et',\et}\et'
\mod v\cl_{R'}$ where $c_{\et',\et}\in\ZZ$.}

\subhead 1.11\endsubhead
Recall that $m=1$. Note that for $\o\in\un{\un\YY}^\bul$ we have $d(\o)\le\sha(R)$. One can show 
(see 8.3):

(a) {\it There is a unique $\o_0\in\un{\un\YY}^\bul$ such that $d(\o_0)=\sha(R)$; moreover,
$\BB^{\o_0}$ consists of a single element $\et_0$.}

\subhead 1.12\endsubhead
Recall that $m=1$. Let $\hW$ be the set of isomorphism classes of irreducible representations 
of $W$ (over $\QQ$). 

For any $E\in\hW$ we denote by $b_E$ the smallest integer $n$ such that $E$ appears in the 
$n$-th symmetric power of the reflection representation of $W$. One can show (see 8.3):

(a) {\it There is a unique bijection $\BB@>\si>>\hW$, $\et\m E_\et$ such that (i),(ii),(iii)
below hold.}

(i) {\it If $R=\emp$ then $\BB@>\si>>\hW$ is the unique bijection between two sets with one
element.}

(ii) {\it If $\r,R',W'$ are as in 1.9 and $R'\ne R$, then for any $\et\in\BB,\et'\in\BB_{R'}$, the integer
$c_{\et',\et}$ in 1.10(b) is equal to the multiplicity of $E_{\et'}$ in $E_\et|_{W'}$. (Here 
we assume that the bijection $\BB_{R'}@>\si>>\hW'$, $\et'\m E_{\et'}$ is already established 
for $R'$ instead of $R$ when $R'\ne R$).}

(iii) {\it For any $\et\in\BB$ we have $(\et:\et_0)=cv^{2b}\mod v^{2b+2}\ZZ[v^2]$ where 
$b=b_{E_\et}$ and $c\in\ZZ_{>0}$.}
\nl
For any $\o\in\un{\un\YY}^\bul$
we denote by $\hW^\o$ the subset of $\hW$ corresponding under (a)
to the subset $\BB^\o$ of $\BB$. From (a) we deduce:

(b) {\it We have $\hW=\sqc_{\o\in\un{\un\YY}^\bul}\hW^\o$ and for any $\o\in\un{\un\YY}^\bul$,
$\et\m E_\et$ in (a) restricts to a bijection $\BB^\o@>\si>>\hW^\o$.}

\subhead 1.13\endsubhead
We now drop the assumption that $m<\iy$; in the case where $m=\iy$ we shall use notation and
results in Appendix B but we will omit the symbol $\d$ which we assume to be $1$. Let
$\o\in\un{\un\YY}^\bul$. 
Let $\r$ be an $m$-facet in $\o$. Then the (rigid) $\ZZ$-graded root system $R(\r)_*$
and the element $y_{R(\r)_*}\in\YY$ are defined.

Let $W(\r)$ be the Weyl group of $R(\r)$; we have $W(\r)\sub W$. Let $W(\r)_0$ be the 
subgroup of $W(\r)$ generated by reflections with respect to roots in $R(\r)_0$; this is 
equal to the stabilizer of $\r$ in $W_0$.

We now disregard for a moment the $\ZZ$-grading of $R(\r)$ and view $R(\r)$ with the obvious
$1$-grading; then the $1$-facets relative to $R(\r)$ are defined and we denote by $\ti\r$ the
$1$-facet relative to $R(\r)$ that contains $y_{R(\r)_*}$. Note that $\ti\r$ is $1$-rigid.
Let $(\ti\r)$ be the $W(\r)_0$-orbit of $\ti\r$.

We set $\hW^{[\o]}=\widehat{W(\r)}^{(\ti\r)}$. This is well defined (independent of the choice of $\r$).
We now consider the bijections
$$\BB^\o\lra Z^\o_R@>f_\r\i>>{}^1Z(\r)^{[0]_\r}_{R(\r)}@>f_{\ti\r}>>Z^{(\ti\r)}_{R(\r)}\lra
\BB^{(\ti\r)}_{R(\r)}\lra\widehat{W(\r)}^{(\ti\r)}$$
where:

the first bijection is as in 1.6(e) (for $\r$, $m<\iy$) or as in B.6(c) (for $\r$, $m=\iy$);

the second bijection is as in 1.6(a) (for $\r$, $m<\iy$) or as in B.6(c) (for $\r$, $m=\iy$);

the third bijection is as in 1.6(a) (for $\ti\r$, $m=1$);

the fourth bijection is as in 1.6(e) (for $\ti\r$, $m=1$);

the fifth bijection is as in 1.12(b) with $R$ replaced by $R(\r)$.
\nl
The composition of these bijections is a bijection 
$\BB^\o\lra\widehat{W(\r)}^{(\ti\r)}$. This can be viewed as a canonical bijection 
$$\BB^\o\lra\hW^{[\o]}.\tag a$$ 
Taking disjoint union over all $\o$ we obtain a bijection
$$\BB\lra\sqc_{\o\in\un{\un\YY}^\bul}\hW^{[\o]}.\tag b$$ 

\head 2. Weyl group representations \endhead
\subhead 2.1\endsubhead
Let $W$ be a Weyl group (not necessarily the one in 1.1); 
we denote by $S$ a set of simple reflections for $W$. Let 
$\hW$ be the set of isomorphism classes of irreducible representations (over $\QQ$) of
$W$. For any $E\in\hW$ let $b_E$ be as in 1.12. For any 
$J\sub S$ let $W_J$ be the subgroup of $W$ generated by $J$.

\proclaim{Proposition 2.2}Let $E,E'$ in $\hW$ be such that (i),(ii) below hold.

(i) $b_E=b_{E'}$;

(ii) for any $J\subsetneqq S$, the restrictions $E|_{W_J}$, $E'|_{W_J}$ are isomorphic.
\nl
Then $E=E'$.
\endproclaim
Assume first that $W=W^1\T W^2,S=S^1\sqc S^2$ where $(W^1,S^1),(W^2,S^2)$ are Weyl 
groups such that $S^1\ne\emp,S^2\ne\emp$ and that the result is known when $(W,S)$ is 
replaced by $(W^1,S^1)$ or by $(W^2,S^2)$. We can write $E=E^1\bxt E^2$,
$E'=E'{}^1\bxt E'{}^2$ with $E^1,E'{}^1$ in $\hW^1$ and $E^2,E'{}^2$ in $\hW^2$.
Taking in (ii) $J=S^1$ we see that $(E^1)^{\op\dim E^2}\cong(E'{}^1)^{\op\dim E'{}^2}$ 
as $W^1$-modules, hence $E^1\cong E'{}^1$ as $W^1$-modules. Similarly, $E^2\cong E'{}^2$ as 
$W^2$-modules. It follows that $E=E'$. Thus we are reduced to the case where $W$ is an 
irreducible Weyl group, hence $|S|\ge1$, which we assume in the remainder of the proof.

For any $n\in\NN$ let $P_n$ be the set of sequences $[\l_1\ge\l_2\ge\do]$ of integers $\ge0$ 
with $\l_k=0$ for large $k$ and with $\sum_k\l_k=n$. Let $(P_n)$ be the group of formal 
$\ZZ$-linear combinations of elements in $P_n$; let $\nu_n:(P_n)@>>>\ZZ$ be the function given 
by the sum of coefficients of an element in $(P_n)$. When $n\ge1$ we define 
$f_n:P_n@>>>(P_{n-1})$ by
$$[\l_1\ge\l_2\ge\do]\m
\sum_{i\ge1;\l_i>\l_{i+1}}[\l_1\ge\l_2\ge\do\ge\l_{i-1}\ge\l_i-1\ge\l_{i+1}\ge\do].$$
For $n=0$ we define $(P_{-1})=0$ and $f_0$ to be the $0$-map. For $n\ge1$ let
$BP_n=\sqc_{n',n''\in\NN;n'+n''=n}P_{n'}\T P_{n''}$. Let $(BP_n)$ be the group of formal $\ZZ$-linear
combinations of elements in $BP_n$; we identify $(BP_n)=\op_{n',n''\in\NN;n'+n''=n}(P_{n'})\ot(P_{n''})$
in an obvious way. When $n\ge2$ we define $\ff_n:BP_n@>>>(BP_{n-1})$ by
$$(\l',\l'')\m f_{n'}(\l')\ot\l''+\l'\ot f_{n''}(\l'')$$
where $(\l',\l'')\in P_{n'}\T P_{n''}$. 

Let $\s:BP_n@>>>BP_n$ be the involution $(\l',\l'')\m(\l'',\l')$. This induces an
involution of $(BP_n)$ denoted again by $\s$.

\mpb

If $W$ is of type $A_1$ then $\hW$ consists of two objects, one with $b=0$ and one 
with $b=1$; hence in this case the desired statement follows from the assumption (i).

Next we assume that $W$ is of type $A_{n-1}$ with $n\ge3$. We show that in this case the 
desired statement follows from assumption (ii) where $J$ is such that $W_J$ has type $A_{n-2}$. 
We identify $\hW$ with $P_n$ in the standard way; in particular, $E,E'$ correspond to 
$[\l_1\ge\l_2\ge\do]$, $[\l'_1\ge\l'_2\ge\do]$ in $P_n$ and from (ii) with $J$ as above we 
have $f_n([\l_1\ge\l_2\ge\do])=f_n([\l'_1\ge\l'_2\ge\do])$ in $(P_{n-1})$. It is enough to show 
that $[\l_1\ge\l_2\ge\do]$ can be recovered from $f_n([\l_1\ge\l_2\ge\do])$. Let 
$$c=\nu_{n-1}(f_n([\l_1\ge\l_2\ge\do])).$$ 
Now $f_n([\l_1\ge\l_2\ge\do])$ is a sum of terms $[\mu_1\ge\mu_2\ge\do]$ where for any $i$, 
$\mu_i$ is either $\l_i$ or $\l_i-1$. If $c\ge2$ then $\l_i$ is the maximum of all $\mu_i$ for 
the various terms as above. If $c=1$ then 
$$[\l_1\ge\l_2\ge\do]=[a,a,\do,a,0,\do]$$ 
for some $a>0$ and 
$$[\mu_1\ge\mu_2\ge\do]=[a,a,\do,a,a-1,0,\do]$$ 
(a sequence with $k\ge1$ nonzero terms); if this sequence contains some entry $\ge2$, then 
$\l_i=\mu_i$ for $i\ne k$ and $\l_i=\mu_i+1$ for $i=k$; if this sequence contains only entries 
$1$ and $0$ then $k\ge2$ (since $n\ge3$) and $\l_i=\mu_i$ for $i\ne k+1$ and $\l_i=\mu_i+1$ 
for $i=k+1$. This proves our claim and completes the proof of the proposition for type $A$.

Now we assume that $W$ is of type $B_2$. In this case the desired statement follows easily 
from assumption (ii) (we must use both $J$ with $\sha(J)=1$).

Next we assume that $W$ is of type $B_n$ or $C_n$ with $n\ge3$. We show that in this case the 
desired statement follows from assumption (ii) where $J$ is such that $W_J$ has type $B_{n-1}$ 
or $C_{n-1}$. The proof borrows some arguments of Shoji \cite{\SHO}. We identify $\hW$ with 
$BP_n$ in the standard way; in particular, $E$ corresponds 
to $(\l',\l'')\in P_{n'}\T P_{n''}$ and $E'$ corresponds to $(\mu',\mu'')\in P_{k'}\T P_{k''}$ where 
$n'+n''=k'+k''=n$. Then from (ii) with $J$ as above we have (as in \cite{\SHO}) 
$\ff_n(\l',\l'')=\ff_n(\mu',\mu'')$. If $n',n''$ are both $\ne0$, it follows immediately that 
$n'=k'$, $n''=k''$ and $\l'=\mu'$, $\l''=\mu''$. If $n'=0$ then we have $k'=0$, $n''=k''=n$ and 
$f_n(\l'')=f_n(\mu'')$; using the argument in the proof for type $A$ we deduce that 
$\l''=\mu''$ hence $(\l',\l'')=(\mu',\mu'')$. Similarly, if $n''=0$ we have 
$(\l',\l'')=(\mu',\mu'')$. This completes the proof of the proposition for type $B,C$.

Next we assume that $W$ is of type $D_n$ with $n\ge4$. We can regard $W$ as a subgroup of 
index $2$ in a Weyl group $\un W$ of type $B_n$ with set $\un S$ of simple reflections. More 
precisely, we can find $s\ne s'$ in $\un S$ such that $ss'$ has order $4$, 

$S=\{s_1\in\un S;s_1\ne s'\}\sqc\{s'ss'\}$.
\nl
Let $\un J\sub\un S$ be such that the subgroup 
$\un W_{\un J}$ generated by $\un J$ is
a Weyl group of type $B_{n-1}$. Let $J\sub S$ be such that $s\in J$, $s'ss'\in J$ and 
$W_J$ has type $D_{n-1}$. Note that $W_J$ has index $2$ in $\un W_{\un J}$.

Let $E_1\in\hat{\un W}$, $E'_1\in\hat{\un W}$ be such that $E$ (resp. $E'$) is contained 
in the restriction of $E_1$ (resp. $E'_1$) to $W$. Then $E_1$ (resp. $E'_1$) 
corresponds as above to an element $(\l',\l'')$ (resp. $(\mu',\mu'')$) of $BP_n$. We 
have $\l'\in P_{n'},\l''\in P_{n''}$, $\mu'\in P_{k'},\mu''\in P_{k''}$ where $n'+n''=k'+k''=n$. 
Let $E_2\in\hat{\un W}$, $E'_2\in\hat{\un W}$ be such that $E_2$ (resp. $E'_2$) 
corresponds as above to $(\l'',\l')$ (resp. $(\mu'',\mu')$). 

Assume first that $\l'\ne\l''$, $\mu'\ne\mu''$. Then $E=E_1|_W,E'=E'_1|_W$. We have
$$\ind_{W_J}^{\un W_{\un J}}(E|_{W_J})=(E_1\op E_2)|_{\un W_{\un J}},\qua
\ind_{W_J}^{\un W_{\un J}}(E'|_{W_J})=(E'_1\op E'_2)|_{\un W_{\un J}}.$$
Since $E|_{W_J}\cong E'|_{W_J}$, we have
$$(E_1\op E_2)|_{\un W_{\un J}}\cong (E'_1\op E'_2)|_{\un W_{\un J}},$$
hence
$$\ff_n(\l',\l'')+\ff_n(\l'',\l')=\ff_n(\mu',\mu'')+\ff_n(\mu'',\mu').$$
If $n',n''$ are both $\ne0$, it follows immediately that either
$n'=k'$, $n''=k''$, $\l'=\mu'$, $\l''=\mu''$ or 
$n'=k''$, $n''=k'$, $\l'=\mu''$, $\l''=\mu'$; in both cases we have $E=E'$.
 
If $n'=0$ then we have either $k'=0$, $n''=k''=n$, $f_n(\l'')=f_n(\mu'')$ (hence
$\l''=\mu''$) or $k''=0$, $k'=n'=n$, $f_n(\mu')=f_n(\l')$ (hence $\mu'=\l'$);
thus we have $E=E'$. Similarly, if $n''=0$ we have $E=E'$.

Next we assume that $\l\ne\ti\l$, $\l'=\ti\l'$ hence $E'_1=E'_2$. In this case we have $\ind_W^{\un W}(E')=E'_1$.
We have
$$\ind_{W_J}^{\un W_{\un J}}(E|_{W_J})=(E_1\op E_2)|_{\un W_{\un J}},\qua
\ind_{W_J}^{\un W_{\un J}}(E'|_{W_J})=E'_1|_{\un W_{\un J}}.$$

Since $E|_{W_J}\cong E'|_{W_J}$, we have
$$(E_1\op E_2)|_{\un W_{\un J}}\cong E'_1|_{\un W_{\un J}},$$
hence
$$\ff_n(\l,\ti\l)+\ff_n(\ti\l,\l)=\ff_n(\l',\l').$$
This is impossible. Similarly if $\l=\ti\l$, $\l'\ne\ti\l'$ we have a contradiction.

We now assume that $\l=\ti\l$, $\l'=\ti\l'$ hence $E_1=E_2$, $E'_1=E'_2$. We have
$$\ind_{W_J}^{\un W_{\un J}}(E|_{W_J})=E_1|_{\un W_{\un J}},\qua
\ind_{W_J}^{\un W_{\un J}}(E'|_{W_J})=E'_1|_{\un W_{\un J}}.$$
Since $E|_{W_J}\cong E'|_{W_J}$, we have $E_1|_{\un W_{\un J}}\cong E'_1|_{\un W_{\un J}},$
hence
$$\ff_n(\l,\l)=\ff_n(\l',\l').$$
It follows that $n'=n''=k'=k''=n/2$, $\l=\l'$.
Hence $E_1=E'_1$. Now $E_1|_W$ splits as a direct sum of two nonisomorphic irreducible
$W$-modules $\ce,\ce'$; $E$ is isomorphic to $\ce$ or to $\ce'$; similarly $E'$ is 
isomorphic to $\ce$ or to $\ce'$. Assume that $E\ne E'$; then $\{E,E'\}=\{\ce,\ce'\}$.
Let $J_1=S-\{s\}$; then $W_{J_1}$ has type $A_{n-1}$. 
We can find $J_2\sub J_1$ such that $\ce|_{W_{J_2}}$ contains the sign representation 
of $W_{J_2}$ but $\ce'|_{W_{J_2}}$ does not contain the sign representation 
of $W_{J_2}$. (This can be deduced from  \cite{\ORA, (4.6.2)}.) If $E\ne E'$ then
it follows that $\ce|_{W_{J_2}}\not\cong\ce'|_{W_{J_2}}$ hence
$E|_{W_{J_2}}\not\cong E'|_{W_{J_2}}$; this contradicts the assumption (ii). We
see that $E=E'$. This completes the proof of the proposition for type $D$.

We now assume that $W$ is of type $G_2$. Note that $\dim E=\dim E'$. A one
dimensional representation of $W$ is determined by its restrictions to the two
$W_J$ with $\sha(J)=1$. Thus if $E,E'$ are $1$-dimensional then $E=E'$ follows from the
assumption (ii). If $E,E'$ are $2$-dimensional and nonisomorphic then their 
$b$-function is $1$ for one of them and $2$ for the other, contradicting the
assumption (i). This completes the proof of the proposition for type $G_2$.

In the remaining cases we shall use the induction/restriction tables of Alvis \cite{A}.

Assume that $W$ is of type $F_4$. From the tables in \cite{\ALV} we see that
$E$ is determined by its restriction to $W_J$ of type $C_3$ or $B_3$. Hence using assumption 
(ii) we must have $E=E'$. This completes the proof of the proposition for type $F_4$.

We now assume that $W$ is of type $E_6$. From the tables in \cite{\ALV} we see that
$E$ is determined by its restriction to $W_J$ of type $A_5$. Hence using assumption 
(ii) we must have $E=E'$. This completes the proof of the proposition for type $E_6$.

We now assume that $W$ is of type $E_7$. From the tables in \cite{\ALV} we see that
$E$ is determined by its restriction to $W_J$ of type $E_6$ except when $\dim E=512$.
Hence using the assumption (ii) we must have $E=E'$ provided that $E,E'$ have 
dimension $\ne512$. Assume now that $E,E'$ have dimension $512$ and are nonisomorphic. 
Then their $b$-function is $11$ for one of them and $12$ for the other, contradicting 
assumption (i). This completes the proof of the proposition for type $E_7$.

We now assume that $W$ is of type $E_8$. From the tables in \cite{\ALV} we see that
$E$ is determined by its restriction to $W_J$ of type $E_7$. Hence using assumption 
(ii) we must have $E=E'$. This completes the proof of the proposition for type $E_8$.
The proposition is proved.

\head 3. Cyclically graded Lie algebras \endhead
\subhead 3.1\endsubhead
In the rest of this paper we assume that $\kk$ is an
algebraic closure of the finite prime field $\FF_p$ with $p$ elements where $p$ is a 
large prime number. For any $q$, a power of $p$, we denote by $\FF_q$ the subfield of
$\kk$ with $\sha(\FF_q)=q$.

If $X$ is an
algebraic variety and $K\in\cd(X)$, $n\in\ZZ$ we write (as in \cite{\GRAII, 3.1})
$K[[n/2]]$ instead of $K[n]\ot\bbq(n/2)$.

Let $Y=\Hom(\kk^*,T)$, $X=\Hom(T,\kk^*)$ viewed as abelian groups with operation 
written as addition; let $(,):Y\T X@>>>\ZZ$ be the obvious pairing. This extends 
to a perfect bilinear pairing $(,):\YY\T\XX@>>>\QQ$ 
where $\YY=\QQ\ot Y$, $\XX=\QQ\ot X$ are viewed as $\QQ$-vector spaces. Let $R\sub X$ (resp. 
$\chR\sub Y$) be the set of roots (resp. coroots) of $G$ with respect to $T$ and let 
$\chR\lra R$ be the standard bijection.
Then $(\YY,\XX,(,),\chR\lra R)$ is a root system as in 1.1(a). Let $\tR=R\cup\{0\}$.
For any $\a\in X$ we set $\fg^\a=\{x\in\fg;\Ad(t)x=\a(t)x\qua\frl t\in T\}$; we have
$\fg^0=\ft$, $\dim\fg^\a=1$ if $\a\in R$, $\fg^\a=0$ if $\a\n\tR$, $\fg=\op_{\a\in\tR}\fg^\a$.

Let $N_GT$ be the normalizer of $T$ in $G$. Let $W=N_GT/T$ be the Weyl group; it can be 
identified with $W$ in 1.1. For any $w\in W$ we denote by $\dw$ a representative of $w$ in 
$N_GT$.

\subhead 3.2\endsubhead
We assume that 

(i) if $m=\iy$, we are given a $\ZZ$-grading $R_*=(R_N)_{N\in\ZZ}$ of $R$ as in B.2.

(ii) if $m<\iy$, we are given a $\ZZ/m$-grading $(R_j)_{j\in\ZZ/m}$ of $R$ as in 1.2.
\nl
If $m=\iy$ we set $\tR_0=R_0\cup\{0\}$ and $\tR_N=R_N$ if $N\ne0$. For $N\in\ZZ$ we set 
$\fg_N=\op_{\a\in\tR_N}\fg^\a$. Then $\fg=\op_{N\in\ZZ}\fg_N$ is a $\ZZ$-grading of $\fg$. If 
$m<\iy$, we set $\tR_0=R_0\cup\{0\}$, $\tR_j=R_j$ if $j\in\ZZ/m-\{0\}$. For $j\in\ZZ/m$ we set
$\fg_j=\op_{\a\in\tR_j}\fg^\a$. Then $\fg=\op_{j\in\ZZ/m}\fg_j$ is a $\ZZ/m$-grading of $\fg$.

For $m\le\iy$ we have $\fg_0=\fL G_0$ where $G_0$ is a closed connected reductive subgroup of 
$G$. If $m=\iy$ we have $\Ad(g)\fg_N=\fg_N$ for $g\in G_0$, $N\in\ZZ$. If $m<\iy$ we have 
$\Ad(g)\fg_j=\fg_j$ for $g\in G_0$, $j\in\ZZ/m$.

The $m$-facets and $m$-alcoves in $\YY$ are defined as in B.1 (when $m=\iy$) and 1.2 (when 
$m<\iy$).

\subhead 3.3\endsubhead
We assume that we are given $q_0\in\{p,p^2,\do\}$ and a rational $\FF_{q_0}$-structure on 
$G$ with Frobenius map $F_0:G@>>>G$ such that $T$ is defined and split over 
$\FF_{q_0}$ and $F_0(\dw)=\dw$ for any $w\in W$; 
then $G_0$, $\fg_1^{nil}$ inherit $\FF_{q_0}$-structures with Frobenius map $F_0$. 
We shall assume, as we may, that each $G_0$-orbit $\co$ in $\fg_1^{nil}$ is
$F_0$-stable and for each $(\co,\ce)\in\ci$ (see 0.2), we are given an isomorphism
$F_0^*\ce\cong\ce$ which makes $\ce$ pure of weight $0$.
Various other varieties associated to $G$ will be considered with the induced
$\FF_{q_0}$-structure.

Let $\r$ be an $m$-facet in $\YY$. If $m=\iy$, for $N\in\ZZ$ we set 
$\fp^\r_N=\op_{\a\in\tR_N;(y,\a)\ge0}\fg^\a$ where $y\in\r$, so that $\op_N\fp^\r_N$ is a 
parabolic subalgebra of $\fg$, and $\fp^\r_0=\fL P^\r$ where $P^\r$ is a parabolic subgroup of 
$G_0$ containing $T$; for $N\in\ZZ$ we set
$$\fu^\r_N=\op_{\a\in\tR_N;(y,\a)>0}\fg^\a.$$
We have $\fp^\r_1\sub\fg_1^{nil}$.

If $m<\iy$, for $N\in\ZZ$ we set $\fp^\r_N=\op_{\a\in\tR_{\bN};(y,\a)\ge N}\fg^\a$ where 
$y\in\r$, so that $(\fp^\r_N)_{N\in\ZZ}$ is a spiral (see \cite{\LYI, 2.5}), and $\fp^\r_0=\fL P^\r$ 
where $P^\r$ is a parabolic subgroup of $G_0$ containing $T$; for $N\in\ZZ$ we set
$$\fu^\r_N=\op_{\a\in\tR_{\bN};(y,\a)>N}\fg^\a,\qua \fl^\r_N=\op_{\a\in\tR_{\bN};(y,\a)=N}\fg^\a$$
so that $\fp^\r_N=\fu^\r_N\op\fl^\r_N$ for all $N$ and $\fl^\r:=\op_N\fl^\r_N=\fL L^\r$,
$\fl^\r_0=\fL L^\r_0$ where 
$L^\r, L^\r_0$ are connected reductive subgroups of $G$ containing $T$. For $N\in\ZZ$ let 
$\p^\r_N:\fp^\r_N@>>>\fl^\r_N$ be the linear map which is $1$ on $\fl^\r_N$ and is $0$ on 
$\fu^\r_N$. We have $\fp^\r_1\sub\fg_1^{nil}$.

For $m\le\iy$ we have $\Ad(g)\fp^\r_1=\fp^\r_1$ for any $g\in P^\r$. We have a diagram
$$\fl^\r_1@<c<<E'@>b>>E''@>a>>\fg_1^{nil}\tag a$$
where 
$$E'=\{(g,z)\in G_0\T\fg_1;\Ad(g\i)z\in\fp^\r_1\},$$ 
$$E''=\{(gP^\r,z)\in G_0/P^\r\T\fg_1;\Ad(g\i)z\in\fp^\r_1\},$$
$$c(g,z)=\p^\r_1(\Ad(g\i)z), b(g,z)=(gP^\r,z), a(gP^\r,z)=z.$$
\nl
We now assume that $\r$ is an $m$-alcove. Then $P^\r$ is a Borel subgroup of $G_0$ and $E''$ 
is smooth, connected of dimension
$$\dim G_0-\dim\fp^\r_0+\dim\fp^\r_1=\dim\fu^\r_0+\dim\fu^\r_1.$$
We set
$$K^\r=a_!\bbq\in\cd(\fg_1^{nil}),
\tK^\r=K^\r[[(\dim\fu^\r_0+\dim\fu^\r_1)/2]]\in\cd(\fg_1^{nil}).$$
Since $a$ is proper, the decomposition theorem \cite{\BBD} shows that $\tK^\r$ is a direct sum 
of shifts of simple perverse sheaves of the form $\ce^\sha[[\dim\co/2]]$ for various 
$(\co,\ce)\in\ci$. Let $\ci_0=\ci_0(\fg_1)$ be the set of all $(\co,\ce)\in\ci$ such 
that some shift of $\ce^\sha$ is a direct summand of $\tK^\r$ for some $m$-alcove $\r$.

For $(\co,\ce)\in\ci_0,(\tco,\tce)\in\ci_0$ we write $(\tco,\tce)<(\co,\ce)$ if 
$\dim\tco<\dim\co$; we write $(\tco,\tce)\le(\co,\ce)$ if either $(\tco,\tce)<(\co,\ce)$ or
$(\tco,\tce)=(\co,\ce)$. Note that $\le$ is a partial order on $\ci_0$.

Let $\cd_0=\cd_0(\fg_1^{nil})$ be the subcategory of $\cd(\fg_1^{nil})$ consisting of complexes $M$ 
such that for any $j$, any composition factor of ${}^pH^j(M)$ is isomorphic to 
$\ce^\sha[[\dim\co/2]]$ for some $(\co,\ce)\in\ci_0$. Let $\ck_0=\ck_0(\fg_1)$ be the free 
$\ca$-module with basis
$$\{\tt_{\co,\ce};(\co,\ce)\in\ci_0\}.$$
If $M\in\cd_0$ has a given mixed structure relative to the $\FF_{q_0}$-structure of 
$\fg_1^{nil}$, we set   
$$\align&gr(M)=\\&
\sum_{(\co,\ce)\in\ci_0,j\in\ZZ,h\in\ZZ}(-1)^j(\text{mult. of 
$\ce^\sha[[\dim\co/2]]$ in }{}^pH^j(M)_h)v^{-h}\tt_{\co,\ce}\in\ck_0.\endalign$$
Here the subscript $h$ denotes the subquotient of pure weight $h$ of a mixed perverse sheaf. 

If $\g$ is an $m$-alcove then (by Deligne's theorem) $\tK^\g$ is a pure complex of weight $0$ 
(with the mixed structure induced by the obvious mixed structure on $\bbq$ which is pure of 
weight $0$). We set
$$\align&I_\g=gr(\tK^\g)\\&=\sum_{(\co,\ce)\in\ci_0,j\in\ZZ}(\text{mult. of 
$\ce^\sha[[\dim\co/2]]$ in }{}^pH^j(\tK^\g))v^{-j}\tt_{\co,\ce}\in\ck_0.\endalign$$

\subhead 3.4\endsubhead
For $(\co,\ce),(\tco,\tce)$ in $\ci$ we define $P_{\tco,\tce;\co,\ce}\in\NN[v\i]$ as in 0.2.
From \cite{\LYII, 13.7(c)} we have
$$(\co,\ce)\in\ci_0,P_{\tco,\tce;\co,\ce}\ne0
\imp(\tco,\tce)\in\ci_0,(\tco,\tce)\le(\co,\ce).\tag a$$
For $(\co,\ce)\in\ci_0$ we denote by $\un\ce$ the extension of $\ce$ to $\fg_1^{nil}$ by 
$0$ on $\fg_1^{nil}-\co$. We show:

(b) {\it If $(\co,\ce)\in\ci_0$ then $\un\ce\in\cd_0$.}
\nl
We argue by induction on $\dim\co$. If $\dim\co=0$, we have $\un\ce=\ce^\sha$ and the result 
follows. Assume now that $\dim\co>0$. We have a distinguished triangle $(\un\ce,\ce^\sha,M)$ 
where $M\in\cd(\fg_1^{nil})$ is such that for any $j\in\ZZ$, the support of $\ch^j(M)$ is 
contained in $\bco-\co$. Moreover from (a) it follows that for any $G_0$-orbit $\co'$ in 
$\bco-\co$, $\ch^j(M)|_{\co'}$ is a local system with all composition factors of the form 
$\ce'$ with $(\co',\ce')\in\ci_0$. Using the induction hypothesis, we see that 
$\ch^j(M)\in\cd_0$. Since this holds for any $j$, it follows that 
$M\in\cd_0$. Using now the distinguished triangle above we deduce that 
$\un\ce\in\cd_0$. This proves (b).

\mpb

We show:

(c) {\it If $M\in\cd_0$, then for any $j\in\ZZ$ and any $G_0$-orbit $\co$ in 
$\fg_1^{nil}$, any composition factor of $\ch^j(M)|_\co$ is of the form $\ce$ with 
$(\co,\ce)\in\ci_0$.}
\nl
We can assume that $M=\ce'{}^\sha$ where $(\co',\ce')\in\ci_0$. In this case the result 
follows from (a).

\subhead 3.5\endsubhead
Fo any $(\tco,\tce)\in\ci_0$ we define an element $\ss_{\tco,\tce}\in\ck_0$ by 
the equations
$$(-v)^{-\dim\co}\tt_{\co,\ce}=\sum_{(\tco,\tce)\in\ci_0;
(\tco,\tce)\le(\co,\ce)}P_{\tco,\tce;\co,\ce}\ss_{\tco,\tce}\tag a$$     
for any $(\co,\ce)\in\ci_0$. (The definition is by induction on $\dim\tco$ using the 
fact that $P_{\tco,\tce;\co,\ce}=1$ if $(\tco,\tce)=(\co,\ce)$.) Note that
$$\{\ss_{\co,\ce};(\co,\ce)\in\ci_0\}\text{ is an $\ca$-basis of }\ck_0.\tag b$$
If $M\in\cd_0$ has a given mixed structure relative to the $\FF_{q_0}$-structure 
of $\fg_1^{nil}$, we set
$$\align&gr'(M)=\\&\sum_{(\tco,\tce)\in\ci_0,j\in\ZZ,h\in\ZZ}(-1)^j
(\text{mult. of $\tce$ in the local system $(\ch^j(M)|_{\tco})_h$})
v^{-h}\ss_{\tco,\tce}\in\ck_0.\endalign$$
Here the subscript $h$ denotes the subquotient of pure weight $h$ of a mixed
local system on $\tco$. Note that $gr'(M[[r/2]])=(-v)^r gr'(M)$ for $r\in\ZZ$. We show:

(c) {\it For $M$ as above we have $gr(M)=gr'(M)$.}
\nl
We can assume that $M=\ce^\sha$ with $(\co,\ce)\in\ci_0$ and $M$ is viewed as a mixed 
complex of pure weight $0$. Using the purity result \cite{\LYII, 12.2}, we see that
$$gr'(M)=\sum_{(\tco,\tce)\in\ci_0;(\tco,\tce)\le(\co,\ce)}
P_{\tco,\tce;\co,\ce}\ss_{\tco,\tce}$$
that is, $gr'(M)=(-v)^{-\dim\co}\tt_{\co,\ce}=gr(M)$. This proves (c).

\subhead 3.6\endsubhead
For any $(\co,\ce)\in\ci_0$ we set
$$\ti\ss_{\co,\ce}=(-v)^{\dim\co}\ss_{\co,\ce}.$$
Let $i:\co@>>>\fg_1^{nil}$ be the inclusion. Now $i_!\ce[[\dim\co/2]]$ is naturally a mixed 
complex (since $\ce$ is pure of weight $0$) and from the definition we have
$\ti\ss_{\co,\ce}=gr'(i_!\ce[[\dim\co/2]])$. Hence 
$$\ti\ss_{\co,\ce}=gr(i_!\ce[[\dim\co/2]]).$$
Note that
$$\{\ti\ss_{\co,\ce};(\co,\ce)\in\ci_0\}\text{ is an $\ca$-basis of }\ck_0.\tag a$$
From 3.5(a) we have
$$\tt_{\co,\ce}=\sum_{(\tco,\tce)\in\ci_0;
(\tco,\tce)\le(\co,\ce)}(-v)^{\dim\co-\dim\tco}P_{\tco,\tce;\co,\ce}\ti\ss_{\tco,\tce}\tag b$$     
where 
$$(-v)^{\dim\co-\dim\tco}P_{\tco,\tce;\co,\ce}=1\text{ if }(\tco,\tce)=(\co,\ce),$$
$$(-v)^{\dim\co-\dim\tco}P_{\tco,\tce;\co,\ce}\in v\ZZ[v]\text{ if }(\tco,\tce)<(\co,\ce)$$ 
(we use the definition of an intersection cohomology complex).

\subhead 3.7\endsubhead
The following result can be deduced from \cite{\GRAI, 17.3} (for $m=\iy$) and from 
\cite{\LYI, 8.4(a)} (for $m<\iy$):

(a) {\it The elements $I_\g$, where $\g\in\un\YY'$, generate the $\QQ(v)$-vector 
space $\QQ(v)\ot_\ca\ck_0$.}

\subhead 3.8\endsubhead
We define a semilinear involution 
$\bar{}:\QQ(v)\ot_\ca\ck_0@>>>\QQ(v)\ot_\ca\ck_0$ by
$$\ov{\tt_{\co,\ce}}=\tt_{\co,\ce}\tag a$$
for all $(\co,\ce)\in\ci_0$. This involution preserves the $\ca$-submodule 
$\ck_0$. We show:

(b) {\it If $\g\in\un\YY'$ then $\ov{I_\g}=I_\g$.}
\nl
An equivalent statement is as follows: for any $j\in\ZZ$ we have 
${}^pH^j(\tK^\g)\cong{}^pH^{-j}(\tK^\g)$. This follows from Deligne's relative hard Lefschetz 
theorem, see \cite{\BBD, 5.4.10}.

\subhead 3.9\endsubhead
Let $\co$ be a $G_0$-orbit in $\fg_1^{nil}$. By a graded analogue of a theorem of 
Morozov-Jacobson-Kostant (see \cite{\GRAI} for $m=\iy$ and \cite{\LYI, 2.3} for $m<\iy$), we can find
elements $e,h,f$ in $\fg$ such that $h\in\ft$, $e\in\co$ (hence $e\in\fg_1$),
$f\in\fg_{-1}$, $[e,f]=h$; moreover, the $W_0$-orbit of $h$ is uniquely determined.
Now $h$ is the differential of an element $y\in Y$. We can view $y$
as an element of $\YY$. Let $C_\co$ be the $W_0$-orbit of $y$ in $\YY$. This is an invariant of
$\co$ and $\co$ can be reconstructed from $C_\co$.
If $m<\iy$ let $\r_y$ be the $m$-facet in $\YY$ that contains $y/2$.
If $m=\iy$ let $\r_y$ be the $\iy$-facet (see B.1) in $\YY$ that contains $-y_{R_*}+y/2$ (with $y_{R_*}$
as in B.2).

\head 4. Parabolic restriction\endhead
\subhead 4.1\endsubhead
In this section we assume that $m=1$. Let $\r$ be an $\iy$-facet in $\YY$ (as in B.1). We set
$$\fp=\op_{\a\in\tR;(y',\a)\ge0}\fg^\a,\fu=\op_{\a\in R;(y',\a)>0}\fg^\a,
\fl=\op_{\a\in\tR;(y',\a)=0}\fg^\a$$
where $y'\in\r$. Then $\fp=\fL P,\fu=\fL U,\fl=\fL L$ where $P$ is a parabolic subgroup of $G$ 
containing $T$, $U$ is the unipotent radical of $P$ and $L$ is the Levi subgroup of $P$ that 
contains $T$. Let $R'=\{a\in R;\fg^\a\sub\fl\}$. Let $\fp^{nil}=\fg^{nil}\cap\fp$,
$\fl^{nil}=\fg^{nil}\cap\fl$. Let 
$\p:P@>>>L$, $\p^{nil}:\fp^{nil}@>>>\fl^{nil}$, be the obvious maps. We define a functor 
$\Res_\r:\cd(\fg^{nil})@>>>\cd(\fl^{nil})$ by $Res_\r(M)=\p^{nil}_!(M|_{\fp^{nil}})$.

Let $\g\in\un\YY'$. We have
$$\fp^\g_1=\op_{\a\in\tR;(y,\a)\ge1}\fg^\a,\qua\fp^\g_0=\op_{\a\in\tR;(y,\a)\ge0}\fg^\a=\fL B,$$
where $B:=P^\g$ is a Borel subgroup of $G$ containing $T$ and $y\in\g$. As in 3.3,
$$E'':=\{(gB,z)\in G/B\T\fg;\Ad(g\i)z\in\fp^\g_1\}@>a>>\fg_1^{nil},$$
where $a(gB,z)\m z$, is a well defined proper morphism and we have
$K^\g=a_!\bbq\in\cd(\fg^{nil})$,
$$\tK^\g=K^\g[[(\dim\fu^\g_0+\dim\fu^\g_1)/2]]\in\cd(\fg^{nil}).$$
Then $\Res_\r(K^\g)=a'_!(\bbq)$ where 
$$\{(gB,z)\in G/B\T\fp;\Ad(g\i)z\in\fp^\g_1\}@>a'>>\fl^{nil}$$
is given by $a'(gB,z)=\p^{nil}(z)$.

Let $W'=N_LT/T\sub W$. Let $\fe$ be a subset of $W$ such that $W=W'\fe$, 
$\sha \fe=\sha W/\sha W'$. We have $G=\sqc_{\e\in\fe}P\dot\e B$. Let $\e\in\fe$. Let
$$S'_\e:=\{(gB,z)\in(P\dot\e B)/B\T\fp;\Ad(g\i)z\in\fp^\g_1\}@>a'_\e>>\fl^{nil}$$ 
be the restriction of $a'$; we set $K'_\e=a'_{\e!}\bbq\in\cd(\fl^{nil})$. Let
$$S''_\e=\{(h(P\cap\Ad(\dot\e)B),z)\in P/(P\cap\Ad(\dot\e)B)\T\fp;\Ad(h\i)z\in\fq\}
@>a''_\e>>\fl^{nil}$$
where $\fq=\Ad(\dot\e)(\fp^\g_1)=\fp^{\e(\g)}$ and $a''_\e(h(P\cap\Ad(\dot\e)B),z)=\p^{nil}(z)$. We have an 
isomorphism $S''_\e@>\si>>S'_\e$, $(h(P\cap\Ad(\dot\e)B),z)\m(h\dot\e B,z)$. Under this 
isomorphism, $a''_\e$ corresponds to $a'_\e$; hence $K'_\e=a''_{\e!}\bbq$. Let
$${}'S_\e=\{(h(L\cap\Ad(\dot\e)B),z)\in P/(L\cap\Ad(\dot\e)B)\T\fp;\Ad(h\i)z\in\fq\}
@>{}'a_\e>>\fl^{nil}$$
where ${}'a_\e(h(P\cap\Ad(\dot\e)B),z)=\p^{nil}(z)$. The map ${}'S_\e@>>>S''_\e$, 
$(h(L\cap\Ad(\dot\e)B),z)\m(h(P\cap\Ad(\dot e)B),z)$ is an affine space bundle with fibres of 
dimension $\dim(U\cap\Ad(\dot\e)B)=\dim(\fu\cap\fp^{\e(\g)}_0)$. (Note that $\e(\g)\in\un\YY'$ 
hence $\fp^{\e(\g)}_0$ is defined.) We deduce:
$${}'K_\e=K'_\e[[-\dim(\fu\cap\fp^{\e(\g)}_0)]]\tag a$$
where ${}'K_\e={}'a_{\e!}\bbq$.

Now $\e(\g)$ is contained in a unique $1$-alcove $\wt{\e(\g)}$ of $\YY$ defined in terms of 
$R'$ instead of $R$; this $1$-alcove defines a spiral $\fp^{\wt{\e(\g)}}$ of $\fl$ and 
complexes $K^{\wt{\e(\g)}}$,
$$\tK^{\wt{\e(\g)}}=K^{\wt{\e(\g)}}[[(\dim\fu^{\wt{\e(\g)}}_0+\dim\fu^{\wt{\e(\g)}}_1)/2]]$$
on $\fl^{nil}$ in the same way as $\g$ defined the spiral $\fp^\g$ of $\fg$ and 
complexes $K^\g$, $\tK^\g$ on $\fg^{nil}$. We have 
$$\fp^{\wt{\e(\g)}}_1=\op_{\a\in R';(\e(y),\a)\ge1}\fg^\a=\fp^{\e(\g)}_1\cap\fl,$$
$$\fp^{\wt{\e(\g)}}_0=\op_{\a\in R'\cup\{0\};(\e(y),\a)\ge0}\fg^\a=\fL B_\e=
\fp^{\e(\g)}_0\cap\fl=\fl\cap\Ad(\e)\fL B,$$
where $B_\e=L\cap\Ad(\e)B$ is a Borel subgroup of $L$ containing $T$ and 
$K^{\wt{\e(\g)}}=a_{\e!}\bbq$ where
$$S_\e:=\{(uB_\e,z)\in L/B_\e\T\fl;\Ad(u\i)z\in\fp^{\wt{\e(\g)}}_1\}@>a_\e>>\fl^{nil}$$
is given by $a_\e(uB_\e,z)=z$. Since $\ad(u\i)z\in\fl$ and 
$\fp^{\wt{\e(\g)}}_1=\fp^{\e(\g)}_1\cap\fl$, we have
$$S_\e:=\{(uB_\e,z)\in L/B_\e\T\fl;\Ad(u\i)z\in\fq\}.$$
We define $c:{}'S_\e@>>>S_\e$ by $(hB_\e,z)\m(\p(h)B_\e,\p^{nil}(z))$. We show:

(b) {\it $c$ is an affine space bundle with fibres of dimension
$\dim\fu+\dim(\fu\cap\fp^{\e(\g)}_1)$.}
\nl
For $(dB_\e,z)\in S_\e$, the fibre $c\i(dB_\e,z)$ is the set of all $(hB_\e,\tz)\in P/B_\e\T\fp$
such that the image of $hB_\e$ under $P/B_\e@>>>L/B_\e$ is $dB_\e$, $\p^{nil}(\tz)=z$ and 
$\Ad(h\i)\tz\in\fq$. We have $P/B_\e=L/B_\e\T U$ and $\fp=\fl\op\fu$. Hence $c\i(dB_\e,z)$ can 
be identified with 
$$\{(u,z_1)\in U\T\fu;\Ad(u\i)\Ad(d\i)(z+z_1)\in\fq\}.\tag c$$ 
It suffices to show that (c) is an affine space of dimension $\dim\fu+\dim(\fu\cap\fq)$. We set 
$\Ad(d\i)z=d'\in\fl\cap\fq$, $\Ad(u\i)Ad(d\i)z_1=z_2\in\fu$; then (c) becomes
$$\{(u,z_2)\in U\T\fu;\Ad(u\i)d'+z_2\in\fq\}.$$
Using the root decomposition $\fg=\op_{\a\in\tR}\fg^\a$ we see that 
$\fp\cap\fq=(\fl\cap\fq)\op(\fu\cap\fq)$; since $\Ad(u\i)d'+z_2\in\fp\cap\fq$, we have
$\Ad(u\i)d'+z_2=\nu+\mu$ where $\nu\in\fl\cap\fq,\mu\in\fu\cap\fq$ are uniquely determined. 
Setting $z_3=\mu-z_2$ we see that (c) becomes 
$$\{(u,z_3,\nu,\mu)\in U\T\fu\T(\fl\cap\fq)\T(\fu\cap\fq);\Ad(u\i)d'=\nu+z_3\}.$$
We have $\Ad(u\i)d'-d'\in\fu$ (since $d'\in\fl'$). Hence $\nu=d'$ and 
$z_3=\Ad(u\i)d'-d'$. Thus (c) can be identified with $\{(u,\mu)\in U\T(\fu\cap\fq)\}$. This 
proves (b).

Since ${}'a_\e=a_\e c$, from (b) we deduce that 
$${}'K_\e=K^{\wt{\e(\g)}}[[-\dim\fu-\dim(\fu\cap\fp^{\e(\g)}_1)]].$$
Combining this with (a) we deduce
$$K'_\e=K^{\wt{\e(\g)}}[[f_\e]]$$
where 
$$\align&f_\e=-\dim\fu-\dim(\fu\cap\fp^{\e(\g)}_1)+\dim(\fu\cap\fp^{\e(\g)}_0)\\&
=-\sha\{\a\in R;(y',\a)>0\}-\sha\{\a\in R;(y',\a)>0,(\e(y),\a)\ge1\}\\&
+\sha\{\a\in R;(y',\a)>0,(\e(y),\a)\ge0\}\endalign$$
where $y'\in\r,y\in\g$.

\subhead 4.2\endsubhead
For any $\e$, we can view $K^{\wt{\e(\g)}}$ (hence also $K^{\wt{\e(\g)}}[[f_\e]]$) as a pure complex of weight zero, 
by Deligne's theorem applied to the proper map $a_\e$. Using an argument in the proof of 
\cite{\CSI, 3.7} for the partition $G=\sqc_{\e\in\fe}(P\dot\e B)$, we deduce that 
$\Res_\r(K^\g)$ is pure of weight $0$, that
$$\Res_\r(K^\g)\cong\op_{\e\in\fe}K^{\wt{\e(\g)}}[[f_\e]]\text{ in }\cd(\fl^{nil}),$$
that $\Res_\r(\tK^\g)$ is pure of weight $0$ and that
$$\Res_\r(\tK^\g)\cong\op_{\e\in\fe}\tK^{\wt{\e(\g)}}[[f'_\e/2]]
\text{ in }\cd(\fl^{nil})\tag a$$
where
$$f'_\e=2f_\e+\dim\fu^\g_0+\dim\fu^\g_1-\dim\fu^{\wt{\e(\g)}}_0-\dim\fu^{\wt{\e(\g)}}_1.$$
We show:
$$\align&f'_\e=-\sha(\a\in R;(y',\a)>0,(\e(y),\a)>1)+\sha(\a\in R;(y',\a)<0,(\e(y),\a)>1)\\&
-\sha(\a\in R;(y',\a)>0,(\e(y),\a)<0)+\sha(\a\in R;(y',\a)<0,(\e(y),\a)<0)\\&
=-\sum_{\a\in R;(y',\a)\ne0,(\e(y),\a)<0\text{ or }(\e(y),\a)>1}\sg(y',\a).\tag b\endalign$$
We have
$$\align&f'_\e=-2\sha\{\a\in R;(y',\a)>0\}-2\sha\{\a\in R;(y',\a)>0,(\e(y),\a)>1\}\\&
+2\sha\{\a\in R;(y',\a)>0,(\e(y),\a)>0\}+\sha\{\a\in R;(y,\a)>1\}\\&+\sha\{\a\in R;(y,\a)>0\}
-\sha\{\a\in R;(y',\a)=0,(\e(y),\a)>1\}\\&-\sha\{\a\in R;(y',\a)=0,(\e(y),\a)>0\}.\endalign$$
Here we substitute
$$\align&\sha\{\a\in R;(y,\a)>1\}+\sha\{\a\in R;(y,\a)>0\}\\&
=\sha\{\a\in R;(y',\a)\le0,(\e(y),\a)>1\}+\sha\{\a\in R;(y',\a)>0,(\e(y),\a)>1\}\\&
+\sha\{\a\in R;(y',\a)\le0,(\e(y),\a)>0\}+\sha\{\a\in R;(y',\a)>0,(\e(y),\a)>0\}.\endalign$$
We obtain
$$\align&f'_\e=-2\sha\{\a\in R;(y',\a)>0\}-2\sha\{\a\in R;(y',\a)>0,(\e(y),\a)>1\}\\&
+\sha\{\a\in R;(y',\a)\le0,(\e(y),\a)>1\}+\sha\{\a\in R;(y',\a)>0,(\e(y),\a)>1\}\\&
-\sha\{\a\in R;(y',\a)=0,(\e(y),\a)>1\}+2\sha\{\a\in R;(y',\a)>0,(\e(y),\a)>0\}\\&
+\sha\{\a\in R;(y',\a)\le0,(\e(y),\a)>0\}+\sha\{\a\in R;(y',\a)>0,(\e(y),\a)>0\}\\&
-\sha\{\a\in R;(y',\a)=0,(\e(y),\a)>0\},\endalign$$
$$\align&f'_\e=-2\sha\{\a\in R;(y',\a)>0\}-\sha\{\a\in R;(y',\a)>0,(\e(y),\a)>1\}\\&
+\sha\{\a\in R;(y',\a)<0,(\e(y),\a)>1\}+2\sha\{\a\in R;(y',\a)>0,(\e(y),\a)>0\}\\&
+\sha\{\a\in R;(y',\a)<0,(\e(y),\a)>0\}+\sha\{\a\in R;(y',\a)>0,(\e(y),\a)>0\}\\&
=-\sha\{\a\in R;(y',\a)>0,(\e(y),\a)>1\}+\sha\{\a\in R;(y',\a)<0,(\e(y),\a)>1\}\\&
-2\sha\{\a\in R;(y',\a)>0,(\e(y),\a)<0\}\\&
+\sha\{\a\in R;(y',\a)<0,(\e(y),\a)>0\} +\sha\{\a\in R;(y',\a)>0,(\e(y),\a)>0\},\endalign$$
$$\align&f'_\e=-\sha\{\a\in R;(y',\a)>0,(\e(y),\a)>1\}+\sha\{\a\in R;(y',\a)<0,(\e(y),\a)>1\}\\&
-2\sha\{\a\in R;(y',\a)>0,(\e(y),\a)<0\}\\&
+\sha\{\a\in R;(y',\a)>0,(\e(y),\a)<0\} +\sha\{\a\in R;(y',\a)>0,(\e(y),\a)>0\}\\&
=-\sha\{\a\in R;(y',\a)>0,(\e(y),\a)>1\}+\sha\{\a\in R;(y',\a)<0,(\e(y),\a)>1\}\\&
-\sha\{\a\in R;(y',\a)>0,(\e(y),\a)<0\}+\sha\{\a\in R;(y',\a)>0,(\e(y),\a)>0\}\endalign$$
and (b) follows.

We define $\ci_0(\fl_1)$ in the same way as $\ci_0(\fg_1)$ but in terms of $L$ 
instead of $G$. (We have $\fl_1=\fl$.) We show:

(c) {\it If $(\co,\ce)\in\ci_0$, then $\Res_\r\ce^\sha$ is a direct sum of shifts of 
complexes of the form $\ce'{}^\sha$ for various $(\co',\ce')\in\ci_0(\fl_1)$.}
\nl
We can find $\g\in\un\YY'$ and $d\in\ZZ$ such that $\ce^\sha$ is a direct summand of 
$\tK^\g[[d]]$ so that $\Res_\r\ce^\sha$ is a direct summand of $\Res_\r\tK^\r[[d]]$. Using (a) 
we deduce that $\Res_\r\ce^\sha$ is a direct sum of shifts of simple perverse sheaves which 
appear in $\tK^{\wt{\e(\g)}}$ for some $\e\in\fe$. This proves (c).

\subhead 4.3\endsubhead
We define $\cd_0(\fl_1^{nil})$, $\ck_0(\fl_1)$ in terms of $L,\ci_0(\fl_1)$ in the 
same way as $\cd_0$, $\ck_0$ were defined in terms of $G,\ci_0(\fg_1)$. 
From 4.2(b) we see that $\Res_\r$ restricts to a functor $\cd_0@>>>\cd_0(\fl_1^{nil})$ 
denoted again by $\Res_\r$. There is a well defined $\ca$-linear map 
$gr\Res_\r:\ck_0@>>>\ck_0(\fl_1)$ such that the following holds: if 
$(\co,\ce)\in\ci_0$ and $\ce^\sha[[\dim\co/2]]$ is viewed as a pure complex of weight $0$
then $(gr\Res_\r)(\tt_{\co,\ce})=gr(\Res_\r\ce^\sha[[\dim\co/2]])$ where 
$\Res_\r\ce^\sha[[\dim\co/2]]$ is viewed as a mixed complex with the mixed structure 
induced from that of $\ce^\sha$.
            
From the results in 4.2 we see that for any $\g\in\un\YY'$ we have
$$(gr\Res_\r)(I_\g)=\sum_{\e\in\fe}v^{f'_\e}I_{\wt{\e(\g)}}.\tag a$$

\subhead 4.4\endsubhead
Let $(\co,\ce)\in\ci_0$, $(\co',\ce')\in\ci_0(\fl_1)$. Let $d=\dim\co$, $d'=\dim\co'$. We 
view $\ce$ as a pure local system of weight $0$ on $\co$. From \cite{\BBD, 5.1.14} we deduce that 
$(\Res_\r(\un\ce))|_{\co'}$ is mixed of weight $\le0$; hence for any $i\in\ZZ$, 
$\ch^i(Res_\r(\un\ce))|_{\co'}$ is mixed of weight $\le i$. By \cite{\ICC, 1.2}, we have 
$\dim(\co\cap(\p^{nil})\i(\co')\le(d-d')/2$ hence $\ch^i(Res_\r(\un\ce))|_{\co'}=0$ if 
$i>d-d'$ and $\ch^{d-d'}(Res_\r(\un\ce))|_{\co'}$ is pure of weight $d-d'$. We denote by
$\fm_{\ce',\ce}$ the multiplicity of $\ce'$ in the local system 
$\ch^{d-d'}(Res_\r(\un\ce))|_{\co'}$. Let $j:\co'@>>>\fl^{nil}$ be the inclusion.

Let $Q\in\ca$ be the coefficient of $\ss_{\co',\ce'}$ in
$$gr'(Res_\r(\un\ce))=gr(Res_\r(\un\ce))=(gr Res_\r)(\ss_{\co,\ce})\in\ck_0(\fl_1)$$
with $gr,gr'$ defined in terms of $\fl$ instead of $\fg$. We have
$$\align&Q=\text{ coeff. of $\ss_{\co',\ce'}$ in }gr'(j_!j^*Res_\r(\un\ce[[d]]))\\&
=(\text{mult. of $\ce'$ in }j^*(\ch^{d-d'}Res_\r(\un\ce))v^{-(d-d')}+c,\endalign$$
where
$$c=\sum_{i<d-d',h\le i}(-1)^i(\text{mult. of $\ce'$ in }
j^*(\ch^i\Res_\r(\un\ce))_h)v^{-h}\in\sum_{h<d-d'}\ZZ v^{-h}.$$
(We use that $d,d'$ are even.) Thus, $Q=\fm_{\ce',\ce}v^{-d+d'}\mod\sum_{h<d-d'}\ZZ v^{-h}$ and
$$Qv^{d-d'}=\fm_{\ce',\ce}\mod v\ZZ[v].\tag a$$

\subhead 4.5\endsubhead
Let $Q''\in\ca$ be the coefficient of $\ti\ss_{\co',\ce'}$ 
in $(gr Res_\r)(\ti\ss_{\co,\ce})\in\ck_0(\fl_1)$ in the basis 3.6(a) for $\fl$ instead 
of $\fg$. We have $Q''=v^{d-d'}Q$ with notation of 4.4. Hence from 4.4(a) we deduce
$$Q''=\fm_{\ce',\ce}\mod v\ZZ[v].\tag a$$

\head 5. The set $\ci'_0$\endhead
\subhead 5.1\endsubhead
In this section we assume that $m<\iy$.
For any Borel subalgebra $\fb$ of $\fg$ containing $\ft$ we denote by $\fu$ the nilradical of $\fb$
and we consider the proper morphism
$$\{(gB_0,z)\in G_0/B_0\T\fg_1;\Ad(g\i)z\in\fu_1\}@>a>>\fg_1^{nil}$$
where $B_0$ is the Borel subgroup of $G_0$ such that $\fL B_0=\fb\cap\fg_0$, $\fu_1=\fu\cap\fg_1$
and $a(gB_0,z)=z$. Let $K'{}^\fb=a_!\bbq$. 
 Let $\ci'_0$ be the set of all $(\co,\ce)\in\ci$ such 
that $\ce$ is a direct summand of $\ch^i(K'{}^\fb)|_\co$ for some $i\in\ZZ$ 
and some $\fb$ as above.
(In the case where $m=1$ this condition on $(\co,\ce)$ appears in Springer's work \cite{\SPR}.)
In this section we prove the following result.

\proclaim{Proposition 5.2} We have $\ci_0=\ci'_0$.
\endproclaim

\subhead 5.3\endsubhead
Let $(\co,\ce)\in\ci'_0$. Then
$\ce$ is a direct summand of $\ch^i(K'{}^\fb)|_\co$ for some $i\in\ZZ$ 
where $\fb$ is as in 5.1. Let $\fu$ be the nilradical of $\fb$.
We can find $y\in Y$ such that 
$$\fb=\op_{\a\in\tR;(y,\a)\ge0}\fg^\a,\fu=\op_{\a\in\tR;(y,\a)>0}\fg^\a, 
\fu\cap\fg_1=\op_{\a\in R_1;(y,\a)>0}\fg^\a.$$
 Since $(y,\a)\in\ZZ$ for any $\a$, we must have $\fu\cap\fg_1=\op_{\a\in R_1;(y,\a)\ge1}\fg^\a$. 
Let $\g$ be the $m$-facet containing $y$; it is an $m$-alcove. We have 
$$\fp^\g_1=\op_{\a\in R_1;(y,\a)\ge1}\fg^\a=\fu\cap\fg_1.$$
 From the definitions we see that 
$K^\g=K'{}^\fb$. Hence $\ce$ is a direct summand of $\ch^i(K^\g)|_\co$ for some $i\in\ZZ$. Using this 
and \cite{\LYII, 13.7(a)}, we deduce that some shift of $\ce^\sha$ is a direct summand of $K^{\g'}$ for some
$m$-alcove $\g'$
hence $(\co,\ce)\in\ci_0$. Thus,

(a) {\it $\ci'_0\sub\ci_0$.}

\subhead 5.4\endsubhead
Let $\co$ be a $G_0$-orbit in $\fg_1^{nil}$. Let $y\in C_\co$ and let $\r=\r_y$ (see 3.9).
Let $L=L^\r,L_0=L^\r_0$, see 3.3. We note the following results. 

(i) {\it Let $\co'$ be the open $L_0$-orbit in $\fl^\r_1$. Then $\co'\sub\co$. There is a unique open 
$P^\r$-orbit $\co''$ in $\fp^\r_1$. We have $\co'\sub\co''$ hence $\co''\sub\co$.}
\nl
(See \cite{\LYI, 2.9(b),(e)}.)

(ii) {\it The map $\ce\m\ce|_{\co'}$ is a 1-1 correspondence between the set of irreducible 
$G_0$-equivariant local systems on $\co$ (up to isomorphism) and the set of irreducible 
$L_0$-equivariant local systems on $\co'$ (up to isomorphism).}
\nl
(See \cite{\LYI, 2.9(c).}

(iii) {\it If $g\in G_0,x\in\co'$ and $\Ad(g\i)(x)\in\fp^\r_1$ then $g\in P^\r$.}
\nl
(See \cite{\LYI, 2.9(d).}

From \cite{\LYI, 7.1(e)} we see that the bijection in (ii) restricts to a bijection
$$\{\ce;(\co,\ce)\in\ci_0\}@>\si>>\{\ce';(\co',\ce')\in\ci_0(\fl^\r_1)\}.\tag a$$

\subhead 5.5\endsubhead
We preserve the setup of 5.4. Let $x\in\co'$. Let $\ce$ be such that $(\co,\ce)\in\ci_0$ and 
let $\ce'=\ce|_{\co'}$. Let $E$ be the irreducible representation of $Z_{G_0}(x)/Z_{G_0}(x)^0$ 
corresponding to $\ce$.
Let $E'$ be the irreducible representation of $Z_{L_0}(x)/Z_{L_0}(x)^0$ corresponding to $\ce'$.
By \cite{\LYI, 2.9(c)}
 we can identify

$Z_{G_0}(x)/Z_{G_0}(x)^0=Z_{L_0}(x)/Z_{L_0}(x)^0$
\nl
and then $E$ 
becomes $E'$. By 5.4(a) we have $(\co',\ce')\in\ci_0(\fl^\r_1)$. It follows that 
there exists a Borel subalgebra $\fb'$ of $\fl^\r$ that contains $\ft$ such that $E=E'$ appears 
in the natural representation of $Z_{L_0}(x)/Z_{L_0}(x)^0$ in $\op_iH^i_c(\cx',\bbq)$ where
$$\cx'=\{gB'_0\in L_0/B'_0;\Ad(g\i)x\in\fb'\};$$
here $B'_0$ is the Borel subgroup of $L_0$ such that $\fL B'_0=\fb'\cap\fl^\r_0$. 
Using \cite{\GRAI, 21.1} we see that $H^i_c(\cx',\bbq)=0$ for $i$ odd hence 

(a) {\it $E=E'$ appears in the virtual representation of $Z_{L_0}(x)/Z_{L_0}(x)^0$ in }
$$\sum_i(-1)^iH^i_c(\cx',\bbq).$$
\nl
Let $\fb$ be a Borel subalgebra of $G$ such that $\fb'\sub\fb$. Let $B_0$ be the Borel subgroup 
of $G_0$ such that $\fL B_0=\fb\cap\fg_0$. We define 
$e:L_0/B'_0@>>>G_0/B_0$ by $gB'_0\m gB_0$. This is well defined since $L_0\sub G_0$ hence
$B'_0\sub B_0$; moreover, $e$ is an imbedding since $B_0\cap L_0=B'_0$. For $t\in\kk^*$ we have 
$y(t)\in T$; we define a $\kk^*$-action on $G_0/B_0$ by $t:gB_0\m y(t)gB_0$. The fixed point set 
of this action is
$$(G_0/B_0)^{\kk^*}=\{gB_0\in G_0/B_0;\Ad(g\i)h\in\fb\}.$$
Note that the image of $e$ is contained in $(G_0/B_0)^{\kk^*}$. (We use that $h$ is contained in 
the centre of $\fl^\r_0$ hence $\Ad(g\i)h=h$ for $g\in L_0$.) Thus $e$ restricts to an imbedding
$e':L_0/B'_0@>>>(G_0/B_0)^{\kk^*}$. This identifies $L_0/B'_0$ with a connected component of
$(G_0/B_0)^{\kk^*}$. (We use that $L_0$ is the centralizer of $h$ in $G_0$.) Now $e'$ restricts 
to an imbedding $\cx'@>>>\cx$ where
$$\cx=\{gB_0\in(G_0/B_0)^{\kk^*};\Ad(g\i)x\in\fb\}.$$
(Note that $\cx$ is well defined since $\Ad(y(t)\i)x=t^{-2}x$ for $t\in\kk^*$.) This imbedding
identifies $\cx'$ with with $e'(L_0/B'_0)\cap\cx$, intersection in $(G_0/B_0)^{\kk^*}$. (We use 
that $g\in L_0,\Ad(g\i)x\in\fb\imp \Ad(g\i)x\in\fb'$ which follows from $\fl^\r\cap\fb=\fb'$.)
Now $\cx'=e'(L_0/B'_0)\cap\cx$ is the intersection of $\cx$ with a connected component of
$(G_0/B_0)^{\kk^*}$ hence $\cx'$ is a union of connected components of $\cx$.
Using this and (a) we deduce that

(b) {\it $E=E'$ appears in the virtual representation of $Z_{L_0}(x)/Z_{L_0}(x)^0$ in }
$$\sum_i(-1)^iH^i_c(\cx,\bbq).$$
\nl
Let $\ti\cx=\{gB_0\in G_0/B_0;\Ad(g\i)x\in\fb\}$. Now $\cx$ is the fixed point set of the
$\kk^*$-action on $\ti\cx$ (the restriction of the $\kk^*$-action on $G_0/B_0$.) Using 
(b) and the fact that the (equivariant) Euler characteristic is preserved by passage to the fixed 
point set of a $\kk^*$-action we deduce that $E=E'$ appears in the virtual 
$Z_{G_0}(x)/Z_{G_0}(x)^0=Z_{L_0}(x)/Z_{L_0}(x)^0$-module $\sum_i(-1)^iH^i_c(\ti\cx,\bbq)$. Hence,
for some $i$, $E=E'$ appears in the $Z_{G_0}(x)/Z_{G_0}(x)^0=Z_{L_0}(x)/Z_{L_0}(x)^0$-module 
$H^i_c(\ti\cx,\bbq)$. In other words, we have $(\co,\ce)\in\ci'_0$. Thus we have proved:

(c) {\it $\ci_0\sub\ci'_0$.}

\subhead 5.6\endsubhead
Proposition 5.2 follows from 5.3(a) and 5.5(c).

\head 6. Inner product\endhead
\subhead 6.1\endsubhead
In this section we assume that $m<\iy$. Let $\g,\g'$ be two $m$-alcoves in $\YY$. Let 
$q=q_0^s$ where $s\ge1$ and let $F=F_0^s:G@>>>G$, $F=F_0^s:\fg@>>>\fg$. We fix a square root 
$\sqrt{q_0}$ of $q_0$ in $\bbq$. Let $\sqrt{q}=(\sqrt{q_0})^s$.

We associate $\fp^\g_N,\fu^\g_N$ ($N\in\ZZ$) and $B=P^\g$ to $\g$ as in 3.3; we associate 
in a similar way $\fp^{\g'}_N,\fu^{\g'}_N$ ($N\in\ZZ$) and $B'=P^{\g'}$ to $\g'$. We define 
functions 
$$\c_\g:(\fg_1^{nil})^F@>>>\bbq, \c_{\g'}:(\fg_1^{nil})^F@>>>\bbq$$
by
$$\c_\g(x)=\sha\{gB^F\in G_0^F/B^F;\Ad(g\i)x\in\fp^\g_1\},$$
$$\c_{\g'}(x)=\sha\{gB'{}^F\in G_0^F/B^F;\Ad(g\i)x\in\fp^{\g'}_1\}.$$
Let $U$ (resp. $U'$) be the unipotent radical of $B$ (resp. $B'$). Let 
$W_0=N_{G_0}T/T\sub W$; this is the same as $W_0$ in 1.2. Let $y\in\g,y'\in\g'$. We show:
$$\align&\sum_{x\in(\fg_1^{nil})^F}\c_\g(x)\c_{\g'}(x)\\&= 
\sha(G_0^F/T^F)\sum_{w\in W_0}q^{\sha(\a\in R_1,(y,\a)\ge1,(w(y'),\a)\ge1)
-\sha(\a\in R_0,(y,\a)\ge0,(w(y'),\a)\ge0)}.\tag a\endalign$$
Let $A$ be the left hand side of (a). We set 
$$\fq_0=\fp^\g_0,\fq_1=\fp^\g_1, \fq'_0=\fp^{\g'}_0,\fq'_1=\fp^{\g'}_1.$$
We have
$$\align&A=\\&
\sha\{(gB^F,g'B'{}^F,x)\in G_0^F/B^F\T G_0^F/B'{}^F\T\fg_1^F;\Ad(g\i)x\in\fq_1, 
\Ad(g'{}\i)x\in\fq'_1\}\\&
=a\i\sha\{(g,g',x)\in G_0^F\T G_0^F\T\fg_1^F;\Ad(g\i)x\in\fq_1,\Ad(g'{}\i)x\in\fq'_1\}\endalign$$
where $a=\sha(B^F)\sha(B'{}^F)$. Setting $x_1=\Ad(g\i)x,h=g'{}\i g$, we have
$$A=a\i\sha\{(g,h,x_1)\in G_0^F\T G_0^F\T\fq_1^F;\Ad(h)x_1\in\fq'_1\}.$$
We have $G_0^F=\sqc_{w\in W_0}B'{}^F\dw B^F$. Hence
$$\align&A=\sha(G_0^F)a\i\T\\&
\sum_{w\in W_0}\sha\{(b,b',x_1)\in B^F\T B'{}^F\T\fq_1^F;
\Ad(b'\dw b)x_1\in\fq'_1\}\sha(B^F\cap\dw\i B'{}^F\dw)\i.\endalign$$
Setting $\Ad(b)x_1=x_2$, we see that
$$\align&A=\\&\fra{\sha(G^F)}{a}\sum_{w\in W_0}\sha\{(b,b',x_2)\in B^F\T B'{}^F\T\fq_1^F;
\Ad(\dw)x_2\in\fq'_1\}\sha(\dw\i B'{}^F\dw\cap B^F)\i,\endalign$$
$$A=\sha(G_0^F/T^F)\sum_{w\in W_0}\sha\{\fq_1^F\cap\Ad(\dw\i)\fq'_1{}^F\}
\sha(U^F\cap\dw\i U'{}^F\dw)\i$$
and (a) follows (after changing $w$ to $w\i$).

\subhead 6.2\endsubhead
In the setup of 6.1 we set for $x\in(\fg_1^{nil})^F$:
$$\ti\c_\g(x)=(\sqrt{q})^{-\dim\fu^\g_0-\dim\fu^\g_1}\c_\g(x),$$
$$\ti\c_{\g'}(x)=(\sqrt{q})^{-\dim\fu^{\g'}_0-\dim\fu^{\g'}_1}\c_{\g'}(x).$$
We show:
$$\sum_{x\in(\fg_1^{nil})^F}\ti\c_\g(x)\ti\c_{\g'}(x)= 
\sha(G_0^F/T^F)(\sqrt{q})^{-2\sha(R_0)}\sum_{w\in W_0}(\sqrt{q})^{-\t(y,w(y'))}.\tag a$$
where $\t(,)$ is as in 1.2.

Using 6.1(a), we see that it is enough to show that for any $w\in W_0$ we have
$$\align&2\sha\{\a\in R_1,(y,\a)\ge1,(w(y'),\a)\ge1\}-2\sha\{\a\in R_0,(y,\a)\ge0,(w(y'),\a)\ge0\}\\&
-\sha\{\a\in R_0,(y,\a)\ge0\}-\sha\{\a\in R_1,(y,\a)\ge1\}-\sha\{\a\in R_0,(y',\a)\ge0\}\\&
-\sha\{\a\in R_1,(y',\a)\ge1\}+2\sha(R_0)=-\sha\{\a\in R_1;((y,\a)-1)((w(y'),\a)-1)<0\}\\&
+\sha\{\a\in R_0;(y,\a)(w(y'),\a)<0\}.\endalign$$
Using the equalities
$$\sha\{\a\in R_1,(y',\a)\ge1\}=\sha\{\a\in R_1,(w(y'),\a)\ge1\},$$
$$\align&\sha(R_0)-\sha\{\a\in R_0,(y',\a)\ge0\}=\sha\{\a\in R_0,(y',\a)<0\}\\&=
\sha\{\a\in R_0,(y',\a)>0\}=\sha\{\a\in R_0,(w(y'),\a)>0\},\endalign$$
$$\sha(R_0)-\sha\{\a\in R_0,(y,\a)\ge0\}=\sha\{\a\in R_0,(y,\a)<0\}=
\sha\{\a\in R_0,(y,\a)>0\},$$ 
and setting $w(y')=y''$, we see that it is enough to show:
$$\align&2\sha\{\a\in R_1,(y,\a)\ge1,(y'',\a)\ge1\}-2\sha\{\a\in R_0,(y,\a)\ge0,(y'',\a)\ge0\}\\&
+\sha\{\a\in R_0,(y,\a)\ge0\}-\sha\{\a\in R_1,(y,\a)\ge1\}\\&
+\sha\{\a\in R_0,(y'',\a)\ge0\}-\sha\{\a\in R_1,(y'',\a)\ge1\}\\&=
-\sha\{\a\in R_1;((y,\a)-1)((y'',\a)-1)<0\}+\sha\{\a\in R_0;(y,\a)(y'',\a)<0\}.\endalign$$
It is enough to show that for $N\in\{0,1\}$ we have
$$\align&2\sha\{\a\in R_{\bN},(y,\a)\ge N,(y'',\a)\ge N\}-\sha\{\a\in R_{\bN},(y,\a)\ge N\}\\&
-\sha\{\a\in R_{\bN},(y'',\a)\ge N\}=-\sha\{\a\in R_{\bN};((y,\a)-N)((y'',\a)-N)<0\}.\endalign$$
This is immediate since $(y,\a)\ne N,(y',\a)\ne N$ for any $\a\in R_{\bN}$.

\subhead 6.3\endsubhead
For any mixed complex $\cm$ over a point we define
$$gr(\cm)=\sum_{j\in\ZZ,h\in\ZZ}(-1)^j\dim (\ch^j(\cm))_hv^{-h}\in\ca.$$
Here the subscript $h$ denotes the subquotient of pure weight $h$ of a mixed
$\bbq$-vector space. Let $$i:\fg_1^{nil}@>>>\fg_1^{nil}\T\fg_1^{nil}$$ be the diagonal and let
$r:\fg_1^{nil}@>>>\text{point}$ be the obvious map. We define an $\ca$-bilinear pairing 
$(:):\ck_0\T\ck_0@>>>\ca$ by the requirement that if $(\co,\ce)\in\ci_0$, 
$(\co',\ce')\in\ci_0$ and
$$M=\ce^\sha[[\dim\co]],M'=\ce'{}^\sha[[\dim\co']]$$
are regarded as pure complexes of weight $0$ so that $r_!i^*(M\bxt M')$ is a mixed complex over
the point, then $$(\tt_{\co,\ce}:\tt_{\co',\ce'})=gr(r_!i^*(M\bxt M')).$$

From 6.2(a) we deduce by an argument entirely similar to that in the proof of 
\cite{\GRAII, 3.11(b)} that 
$$(I_\g:I_{\g'})=e_{W_0}\sum_{w\in W_0}v^{\t(y,w(y'))},\tag a$$
with $e_{W_0}$ as in 1.2. Alternatively, one can prove (a) using arguments in the proof of \cite{\LYI, 6.4}.

\subhead 6.4\endsubhead
Let $(\co,\ce)\in\ci_0$, $(\co',\ce')\in\ci_0$. We regard $\un\ce,\un\ce'$ as 
mixed complexes such that $\un\ce|_\co$ and $\un\ce'|_{\co'}$ are pure of weight $0$; then 
$\un\ce\ot\un\ce'$ is a mixed complex and from the definitions we have
$$(\ss_{\co,\ce}:\ss_{\co',\ce'})=gr(r_!(\un\ce\ot\un\ce')).$$
Hence if $\co\ne\co'$ we have
$$(\ss_{\co,\ce}:\ss_{\co',\ce'})=0,\tag a$$
while if $\co=\co'$ we have
$$(\ss_{\co,\ce}:\ss_{\co,\ce'})
=\sum_{j\in\ZZ,h\in\ZZ}(-1)^j\dim((H^j_c(\co,\ce\ot\ce'))_h)v^{-h}\in\ZZ[v\i].\tag b$$
Here the subscript $h$ denotes the subquotient of pure weight $h$ of a mixed
$\bbq$-vector space. Let $d=\dim\co$ and let $\ce^*$ be the local system dual to $\ce$.
Let $\d_{\ce^*,\ce'}$ be $1$ if $\ce'=\ce^*$ and $0$ if $\ce'\ne\ce^*$. We have
$$(H^j_c(\co,\ce\ot\ce'))_h\ne0\imp h\le j\le2d;$$
moreover,
$$(H^{2d}_c(\co,\ce\ot\ce'))_h\ne0\imp h=2d,\d_{\ce^*,\ce'}=1=\dim H^{2d}_c(\co,\ce\ot\ce').$$
It follows that
$$(\ss_{\co,\ce}:\ss_{\co,\ce'})=\d_{\ce^*,\ce'}v^{-2d}\mod v^{-2d+1}\ZZ[v],\tag c$$
so that
$$(\ti\ss_{\co,\ce}:\ti\ss_{\co,\ce'})=\d_{\ce^*,\ce'}\mod v\ZZ[v].$$
From (a),(c) we see that the square matrix $(\ss_{\co,\ce}:\ss_{\co',\ce'})$
has nonzero determinant  hence is invertible over $\QQ(v)$. We deduce that

(d) {\it the $\QQ(v)$-bilinear form
$(:):(\QQ(v)\ot_\ca\ck_0)\T(\QQ(v)\ot_\ca\ck_0)@>>>\QQ(v)$
deduced from $(:):\ck_0\T\ck_0@>>>\ca$ by extension of scalars is
nonsingular.}

\subhead 6.5\endsubhead
Let $i_1=(\co,\ce)\in\ci_0$, $i_2=(\co',\ce')\in\ci_0$. From 3.5(a) we deduce
$$(-v)^{-\dim\co-\dim\co'}(\tt_{i_1}:\tt_{i_2})=\sum_{i'_1\in\ci_0;i'_2\in\ci_0;
i'_1\le i_1,i'_2\le i_2}P_{i'_1,i_1}P_{i'_2,i_2}(\ss_{i'_1}:\ss_{i'_2}).\tag a$$

\subhead 6.6\endsubhead
In the remainder of this section we assume that $m=1$. For $i,i'$ in $\ci_0$ let 
$\Om_{i,i'}\in\QQ[\qq]$ be as in \cite{\CSV, p.145} and let $\L_{i,i'}\in\QQ[\qq]$, 
$\Pi_{i,i'}\in\QQ[\qq]$ be as in \cite{\CSV, p.146}. (Here $\qq$ is an indeterminate.) We shall 
regard $\Om_{i,i'},\L_{i,i'},\Pi_{i,i'}$ as elements of $\QQ(v)$ via $\qq=v^{-2}$.

If $v$ is specialized to $-\sqrt{q}\i$ where $q=q_0^s$, $s\ge1$, then $\qq$ becomes $q^s$,
$\L_{i,i'}$ becomes an integer $\l_{i,i'}$ (depending on $s$) as in the proof in 
\cite{\CSV, p.146}. From the definition, that integer is equal to the specialization of 
$(\ss_{i'_1}:\ss_{i'_2})\in\ZZ[v\i]$ at $v=\sqrt{q}\i$ (we use Grothendieck's trace formula to 
evaluate 6.4(b) and we use that the relevant eigenvalues of Frobenius are integer powers of 
$q$). It follows that
$$(\ss_i:\ss_i')=\L_{i,i'}.\tag a$$
Moreover from the definitions we have
$$P_{i,i'}=\Pi_{i,i'}\tag b$$
where $P_{i,i'}$ is as in 3.4.

Now let $i_1=(\co,\ce)\in\ci_0$, $i_2=(\co',\ce')\in\ci_0$. Using (a),(b), from 6.5(a)
we deduce:
$$v^{-\dim\co-\dim\co'}(\tt_{i_1}:\tt_{i_2})=
\sum_{i'_1\in\ci_0;i'_2\in\ci_0;i'_1\le i_1,i'_2\le i_2}
\Pi_{i'_1,i_1}\Pi_{i'_2,i_2}\L_{i'_1,i'_2}=\Om_{i_1,i_2},\tag c$$
where the last equality follows from equation (b) in \cite{\CSV, p.146}; we have used that
$\dim\co,\dim\co'$ are even.

\subhead 6.7\endsubhead
Let $\hW$ be as in 1.12. Let $i\m E_i$, $\ci'_0@>\si>>\hW$ be the Springer correspondence 
(we use the normalization in \cite{\ICC}). Using the equality $\ci'_0=\ci_0$ in 5.2 we 
can view this as a bijection $\ci_0@>\si>>\hW$. Let $E_0\in\hW$ be the reflection 
representation and let $\nu$ be the dimension of the flag manifold of $G$. Let $e_W=e_{W_0}\in\ca$ 
be as in 1.2. From the definition in \cite{\CSV, 24.7}, for $i,i'$ in $\ci_0$ we have
$$\align&\Om_{i,i'}=\sha(W)\i\sum_{w\in W}\tr(w,E_i)\tr(w,E_{i'})\\&
(v^{-2}-1)^{\dim E_0}\det(v^{-2}-w,E_0)\i\bb(e_W)v^{-\dim\co-\dim\co'+2\nu}.\tag a\endalign$$
Combining this with 6.6(c) we obtain
$$(\tt_i:\tt_{i'})=\sha(W)\i\sum_{w\in W}\tr(w, E_i\ot E_{i'})
(1-v^2)^{\dim E_0}\det(1-v^2w,E_0)\i e_W.\tag b$$
Let $i_0\in\ci_0$ be the element such that $E_{i_0}$ is the unit representation of $W$. 
It is known that $i_0=(\co',\bbq)$ where $\co'$ is the regular nilpotent orbit.
Hence for $i'=i_0$, (b) becomes 
$$(\tt_i:\tt_{i_0})=\sha(W)\i\sum_{w\in W}\tr(w,E_i)(1-v^2)^{\dim E_0}\det(1-v^2w,E_0)\i e_W.
\tag c$$
The right hand side of (c) is the {\it fake degree} $FD(E_i)$ of $E_i$ (see \cite{\MADI, 3.17}). 
Thus we have
$$(\tt_i:\tt_{i_0})=FD(E_i).\tag d$$
We have
$$FD(E_i)=cv^{2b}\mod v^{2b+2}\ZZ[v^2]\tag e$$
where $b=b_{E_i}\in\NN$ is as in 1.12 and $c\in\ZZ_{>0}$ is well defined. From (d) we deduce
$$(\tt_i:\tt_{i_0})=cv^{2b}\mod v^{2b+2}\ZZ[v^2]\tag f$$
where $b=b_{E_i}$ and $c\in\ZZ_{>0}$.

\subhead 6.8\endsubhead
We show:

(a) {\it There exists a unique bijection $\ci_0@>\si>>\hW$, $i\m E'_i$, such that
(i),(ii),(iii) below hold.}

(i) {\it If $R=\emp$ then $\ci_0@>\si>>\hW$ is the unique bijection between two sets 
with one element.}

(ii) {\it Let $\r,\fp,\fl,P,L,W',R'$ be as in 4.1 with $R'\ne R$, and let $i=(\co,\ce)\in\ci_0$, 
$i'=(\co',\ce')\in\ci_0(\fl_1)$. Let $Q''\in\ca$ be as in 4.5. Let $\ti\fm_{i,i'}\in\NN$ be
the multiplicity of $E'_{i'}$ in $E'_i|_{W'}$. Then $Q''=\ti\fm_{i',i}\mod v\ZZ[v]$. (Here 
we assume that the bijection $\ci_0(\fl_1)@>\si>>\hW'$, $i'\m E'_{i'}$ is already established 
for $R'$ instead of $R$ when $R'\ne R$).}

(iii) {\it For any $i\in\ci_0$ we have $(\tt_i:\tt_{i_0})=cv^{2b}\mod v^{2b+2}\ZZ[v^2]$
where $b=b_{E'_i}$ and $c\in\ZZ_{>0}$.}
\nl
If we take $E'_i=E_i$ (see 6.7) then (i) is obvious, (ii) follows from 4.5(a) together with
\cite{\ICC, 8.3(b)} and (iii) follows from 6.7(f). Thus, a bijection as in 
(a) exists. The uniqueness of a bijection as in (a) follows from 2.2. Thus, (a) holds.

\head 7. Induction\endhead
\subhead 7.1\endsubhead
In this section we assume that $m<\iy$. Let $\r$ be an $m$-facet. Let 

$\fp^\r_N,\fu^\r_N,\fl^\r_N,\fl^\r$, $L^\r,L^\r_0,P^\r,E',E''$, $a,b,c$ 
\nl
be as in 3.3. Let $R(\r)$ be the 
set of roots of $L^\r$ with respect to $T$.
This is the same as $R(\r)$ in 1.3. It has a $\ZZ$-grading as in B.4 and the corresponding $\ZZ$-grading
of $\fl^\r$ is given by $\op_N\fl^\r_N$. 

For any $\iy$-alcove $\g$ of $\YY$ with respect to $R(\r)$ we can consider the parabolic 
subalgebra $\op_{N\in\ZZ}\fp^\g_N$ of $\fl^\r$ defined as in 3.3 with $\fg$ replaced by 
$\fl^\r$. (We have $\fp^\g_N\sub\fl^\r_N$ for any $N$.) There is a well defined $m$-alcove 
$\ti\g$ of $\YY$ such that $\fp^{\ti\g}_N=\fp^\g_N\op\fu^\r_N$ for any $N$. This follows from 
the analysis in \cite{\LYI, 2.8} which shows also that $\ti\g=f_\r(\g)$ with $f_\r$ as in 1.4.
Now the complex $\tK^{\g}$ on $\fl^\r_1$ 
(analogous to $\tK^\r$ in 3.3) is defined in terms of the $\ZZ$-grading of $\fl^\r$. Similarly, 
the group $\ck_0(\fl^\r_1)$ is defined in terms of this $\ZZ$-grading and its elements 
$I_\g=gr(\tK^\g)\in\ck_0(\fl^\r_1)$ are defined for any $\g$ as above. Moreover, 
$\ck_0(\fl^\r_1)$ has an $\ca$-basis 
$$\{\tt_{\co',\ce'};(\co',\ce')\in\ci_0(\fl^\r_1)\}$$ 
and an $\ca$-basis 
$$\{\ti\ss_{\co',\ce'};(\co',\ce')\in\ci_0(\fl^\r_1)\}$$
defined as in 3.3, 3.5 (with $\fg$ replaced by $\fl$ with its $\ZZ$-grading).
Here $\ci_0(\fl^\r_1)$ is the set of pairs consisting of an $L^\r_0$-orbit $\co'$ on $\fl^\r_1$ 
and an irreducible $L^\r_0$-equivariant local system $\ce'$ on $\co'$ (up to isomorphism), 
defined like $\ci_0$ (in the $\ZZ$-graded case) but with $G$ replaced by $L^\r$.

Let $\fA$ be a direct sum of shifts of $L^\r_0$-equivariant simple perverse sheaves on 
$\fl^\r_1$. In the diagram 3.3(a), $c$ is smooth with connected fibres of dimension
$\dim G_0+\dim\fu^\r_1$, $b$ is a principal $P^\r$-bundle. Hence there is a well defined (up to 
isomorphism) complex $\fA''$ on $E''$ which is a direct sum of shifts of simple perverse sheaves 
such that
$$c^*\fA[[(\dim G_0+\dim\fu^\r_1)/2]]\cong b^*\fA''[[\dim P^\r/2]].$$
We set $\ti\fA=a_!\fA''$. Since $a$ is proper, $\ti\fA$ is a direct sum of shifts of 
(necessarily $G_0$-equivariant) simple perverse sheaves on $\fg^{nil}_1$. By \cite{\LYI, 4.2}, 
if $\g$ is an $\iy$-facet of $\YY$ with respect to $R(\r)$ and $\ti\g$ is the $m$-facet in 
$\YY$ defined by $\ti\g=f_\r(\g)$ (see 1.4), then applying the previous construction to 
$\fA=\tK^\g$ gives us $\fA''=\tK^{\ti\g}$. As in \cite{\GRAII, 3.5}, it follows that there is a 
well defined $\ca$-linear map $i_\r:\ck_0(\fl^\r_1)@>>>\ck_0$ such that the following 
holds: if $(\co',\cl')\in\ci_0(\fl^\r_1)$ and $\fA=\cl'{}^\sha[[\dim\co'/2]]$, viewed as a pure 
complex of weight $0$, then $\ti\fA$ is canonically defined and
$i_\r(\tt_{\co',\ce'})=gr(\ti\fA)$ where $\ti\fA$ is viewed as a pure complex of weight $0$ 
with mixed structure induced by that of $\fA$. Moreover, if $\g,\ti\g=f_\r(\g)$ are as above, we 
have 
$$i_\r(I_\g)=I_{\ti\g}.\tag a$$

\subhead 7.2\endsubhead
In this subsection we fix a $G_0$-orbit $\co$ in $\fg_1^{nil}$. 
Let $y\in C_\co,\r,\co',\co''$ be as in 5.4. Let 
$$E''_1=\{(gP^\r,z)\in G_0/P^\r\T\fg_1;\Ad(g\i)z\in\co''\},$$
an open subvariety of $E''$ in 3.3(a). We show:

(a) {\it The map $E''_1@>>>\co$, $(gP^\r,z)\m z$ is a well defined isomorphism.}
\nl
This map is well defined since, by 5.4(i), if $(gP^\r,z)\in E''_1$, then $z\in\co$. We shall only 
prove that our map is bijective. Let $z\in\co$. Since $\co''\sub\co$ (by 5.4(i)) and the 
$G_0$-action on $\co$ is transitive, we have $\Ad(g)\i z\in\co''$ for some $g\in G_0$; this 
proves surjectivity of our map. Assume now that $(gP^\r,z)\in E''_1$, $(g'P^\r,z)\in E''_1$. 
Setting $g'=gg_1$ with $g_1\in G_0$ and $z'=\Ad(g\i)z$, we have $z'\in\co''$, 
$\Ad(g_1\i)z'\in\co''$. From 5.4(i) 
we have $z'=\Ad(g_2)x$ where $g_2\in P^\r$, $x\in\co'$. We have 
$\Ad(g_1\i g_2)x\in\co''$; hence using 5.4(iii)
 we have $g_2\i g_1\in P^\r$ that is, $g_1\in P^\r$. 
Thus, $gP^\r=g'P^\r$. This proves that our map is injective hence bijective.

\mpb

Let $\ce$ be such that $(\co,\ce)\in\ci_0$ and let $\ce'=\ce|_{\co'}$ so that 
$(\co',\ce')\in\ci_0(\fl^\r_1)$, see 5.4(a). Then the element $\ti\ss_{\co,\ce}\in\ck_0$ is 
well defined and the analogously defined element $\ti\ss_{\co',\ce'}\in\ck_0(\fl^\r_1)$ is well 
defined. We have
$$i_\r(\ti\ss_{\co',\ce'})=\ti\ss_{\co,\ce}.\tag b$$
The proof is entirely similar to that of \cite{\GRAII, 3.15(d)} (a $\ZZ$-graded analogue
of (b)), using (a) instead of \cite{\GRAII, 3.15(c)}.

\head 8. Proofs\endhead
\subhead 8.1\endsubhead
In this section we finish the proofs of the theorems stated in the Introduction.
We can identify $\QQ(v)\ot_\ca\ck_0$ with $\VV$ (in B.3, if $m=\iy$ or in 1.2 if $m<\iy$)
in such a way that for any $m$-alcove $\g$, the element $I_\g$ of 
$\QQ(v)\ot_\ca\ck_0$ corresponds to the element $I_\g$ of $\VV$ (if $m=\iy$, see B.3) or to 
the element $I_\g$ of $\VV$ (if $m<\iy$, see 1.2). (If $m=\iy$, this follows from the results in 
\cite{\GRAII}. If $m<\iy$ this follows from 3.7(a), 6.3(a), 6.4(d).) Then $(:)$ on
$\QQ(v)\ot_\ca\ck_0$ corresponds to $(:)$ on $\VV$.
Moreover, $\bar{}:\QQ(v)\ot_\ca\ck_0@>>>\QQ(v)\ot_\ca\ck_0$ in 3.8 corresponds to 
$\b:\VV@>>>\VV$ (as in B.3 if $m=\iy$ or as in 1.2(a) if $m<\iy$). This verifies the assertion in
1.2(a).

\subhead 8.2\endsubhead 
Now assume that $m=\iy$. The statements in this subsection
are proved in \cite{\GRAII}. Under the identification in 8.1, the basis 
$\{\ti\ss_{\co,\ce};(\co,\ce)\in\ci_0\}$ of $\QQ(v)\ot_\ca\ck_0$ in 3.6 corresponds to 
the basis ${}^1Z_R$ (see B.6) of $\VV$  and the basis 
$\{\tt_{\co,\ce};(\co,\ce)\in\ci_0\}$ of $\QQ(v)\ot_\ca\ck_0$ in 3.3
corresponds to the basis $\BB$ (see B.6) of $\VV$.

From the results in \cite{\GRAII} we see that the following holds.

(a) {\it Let $\co$ be the open $G_0$-orbit in $\fg_1$ and let $y\in C_\co$, see 3.9. Then
$R_*$ is rigid (see B.3) if and only if $y/2=y_{R_*}$ (with $y_{R_*}$ as in B.2).}
\nl
(Note that the condition that $y/2=y_{R_*}$ is independent of the choice of $y$ in $C_\co$ since
$y_{R_*}$ is fixed by $W_0$.)

It follows that, by associating to $R_*$ the $G$-orbit of $e$ in (a), we get a well defined 
bijection between the set of rigid $\ZZ$-gradings of $R$ (up to $W$-action) and the set of even
nilpotent $G$-orbits in $\fg$.

Assume now that $R_*$ is rigid. Let $\co'$ be the $G$-orbit in $\fg$ corresponding to $R_*$.
Let $\co$ be the open $G_0$-orbit in $\fg_1$, so that $\co\sub\co'$. Then the subset
$\{\ti\ss_{\co,\ce};\ce \text{ such that }(\co,\ce)\in\ci_0\}$ of $\QQ(v)\ot_\ca\ck_0$ 
corresponds to the subset ${}^1Z^{[0]}_R$ (see B.6) of $\VV$.

We show:

(b) {\it An $\iy$-facet $\r$ is $1$-rigid if and only if for some $G_0$-orbit $\co$ on $\fg_1$ we have
$\r=\r_y$ where $y\in C_\co$. Moreover $\co\m\r_y$ is a bijection between the set of
$G_0$-orbits on $\fg_1$ and the set of $W_0$-orbits of $1$-rigid $\iy$-facets in $\YY$.}
\nl
Assume that $\r$ is a $1$-rigid $\iy$-facet (see B.7). Then $y_{R(\r)_*}-y_{R_*}\in\r$ and $R(\r)_*$
is rigid (here $R(\r)_*$ is as in B.4).
Let $\fg(\r)=\op_{\a\in R(\r)\cup\{0\}}\fg^\a$, $\fg(\r)_0=\op_{\a\in R(\r)_0\cup\{0\}}\fg^\a$,
$\fg(\r)_N=\op_{\a\in R(\r)_N}\fg^\a$ if $N\in\ZZ-\{0\}$. Let $G(\r)_0$ be the closed connected subgroup
of $G$ such that $\fL G(\r)_0=\fg(\r)_0$. Then $G(\r)_0$ acts on $\fg(\r)_1$ by $\Ad$. Let $\co_0$
be the open orbit for this action and let $y\in C_{\co_0}$ (defined as in 3.9 in terms of
$\fg(\r)$ instead of $\fg$). Let $\co$ be the $G_0$-orbit on $\fg_1$ that contains $\co_0$.
Note that we have also $y\in C_\co$ (defined as in 3.9 in terms of $\fg$).
By (a) we have $y_{R(\r)_*}=y/2$. Since $y_{R(\r)_*}-y_{R_*}\in\r$ it follows that $y/2-y_{R_*}\in\r$
hence $\r=\r_y$ (see 3.9).
 
Conversely,
assume that $\co$ is a $G_0$-orbit on $\fg_1$. We associate $e,h,f,y,\r_y$ to $\co$ as in 3.9.
Let $\r=\r_y$. Let $R(\r)_*$ be as in B.4. 
Since $y/2-y_{R_*}\in\r$, for $N\in\ZZ$ we have
$R(\r)_N=\{\a\in R_N;(y/2-y_{R_*},\a)=0\}=\{\a\in R_N;(y/2,\a)=N\}$.
We define $\fg(\r),\fg(\r)_N,G(\r)_0$ as above. We have
$\fg(\r)_N=\{x\in\fg_N;[h/2,x]=Nx\}$. Hence $e\in\fg(\r)_1$,
$f\in\fg(\r)_{-1}$. It follows that $e$ is in the open $G(\r)_0$-orbit on $\fg(\r)_1$
and the element $y$ associated as in 3.9 to this open orbit is the same as $y$ above. From the
definitions we have $y/2=y_{R(\r)_*}$. Using (a) we deduce that
$R(\r)_*$ is rigid. We have $y_{R(\r)_*}-y_{R_*}=y/2-y_{R_*}\in\r$. We see that $\r$ is
$1$-rigid. Now (b) follows.

\subhead 8.3\endsubhead 
We now assume that $m<\iy$. Now 1.5(a) follows from 7.1(a); 1.6(a) follows from 7.2(b). We show:

(a) {\it An $m$-facet $\r$ is rigid if and only if for some $G_0$-orbit $\co$ on $\fg_1^{nil}$ we have
$\r=\r_y$ where $y\in C_\co$. Moreover $\co\m\r_y$ is a bijection between the set of
$G_0$-orbits on $\fg_1^{nil}$ and the set of $W_0$-orbits of rigid $m$-facets in $\YY$.}
\nl
The proof is almost a repetition of that of 8.2(b).
Assume that $\r$ is a rigid $m$-facet. Then $y_{R(\r)_*}\in\r$ and $R(\r)_*$ is rigid
in the sense of B.7.
Let $\fg(\r)=\op_{\a\in R(\r)\cup\{0\}}\fg^\a$, $\fg(\r)_0=\op_{\a\in R(\r)_0\cup\{0\}}\fg^\a$,
$\fg(\r)_N=\op_{\a\in R(\r)_N}\fg^\a$ if $N\in\ZZ-\{0\}$. Let $G(\r)_0$ be the closed
connected subgroup of $G$
such that $\fL G(\r)_0=\fg(\r)_0$. Then $G(\r)_0$ acts on $\fg(\r)_1$ by $\Ad$. Let $\co_0$
be the open orbit for this action and let $y\in C_{\co_0}$ (defined as in 3.9 in terms of
$\fg(\r)$ instead of $\fg$). Let $\co$ be the $G_0$-orbit on $\fg_1^{nil}$ that contains $\co_0$.
Note that we have also $y\in C_\co$ (defined as in 3.9 in terms of $\fg$).
By 8.2(a) we have $y_{R(\r)_*}=y/2$. Since $y_{R(\r)_*}\in\r$ it follows that $y/2\in\r$
hence $\r=\r_y$ (see 3.9).

Conversely, assume that $\co$ is a $G^0$-orbit on $\fg_1^{nil}$.
We associate $e,h,f,y,\r_y$ to $\co$ as in 3.9.
Let $\r=\r_y$. Let $R(\r)_*$ be as in 1.3. Since $y/2\in\r$,
for $N\in\ZZ$ we have $R(\r)_N=\{a\in R_{\bN};(y/2,\a)=N\}$.
We define $\fg(\r),\fg(\r)_N,G(\r)_0$ as above. We have
$\fg(\r)_N=\{x\in\fg_{\bN};[h/2,x]=Nx\}$. Hence $e\in\fg(\r)_1$, $f\in\fg(\r)_{-1}$.
It follows that $e$ is in the open $G(\r)_0$-orbit on $\fg(\r)_1$ and the element $y$
associated as in 3.9
to this open orbit is the same as $y$ above. From the definitions we have
$y/2=y_{R(\r)_*}$. Using 8.2(a) we deduce that $R(\r)_*$ is rigid. We have $y_{R(\r)_*}=y/2\in\r$.
We see that $\r$ is rigid.
Now (a) follows.

Let $\co$ be a $G_0$-orbit in $\fg_1^{nil}$. Let $\o$ be the $W_0$-orbit on the set of rigid
$m$-facets associated to $\co$ in (a); then the
subset ${}^1Z^\o_R$ of $\VV$ is defined (see 1.6). Using 7.2(b) we deduce:

(b) {\it Under the identification $\VV=\QQ(v)\ot_\ca\ck_0$, ${}^1Z^\o={}^1Z^\o_R$ becomes
the subset 
$$\{\ti\ss_{\co,\ce};\ce \text{ such that }(\co,\ce)\in\ci_0\}$$ of $\QQ(v)\ot_\ca\ck_0$.}
\nl
Now 1.6(b) follows immediately from (b) and 1.6(c) follows from (a) and 3.6(a).
We show that 1.6(d) holds for any $\x\in Z$. Now $\x$ is of the form $\ti\ss_{\co,\ce}$ for
some $(\co,\ce)\in\ci_0$. Then $\un\x$ in 1.6(d) exists: we can take $\un\x=\tt_{\co,\ce}$
(we use 3.6(b) and 3.8(b)). The uniqueness in 1.6(d) is immediate. Now 1.6(e) follows.

\mpb

Let $\co,\o$ be as above. We show: 

(c) {\it $d(\o)$ in 1.6 is equal to $\dim\co$.}
\nl
By 7.2(a) it is enough to show that $d(\o)$ is equal to $\dim E''_1$ (in 7.2 with $\r\in\o$)
that is to

$\dim(G_0)-\dim\fp^\r_0+\dim\fp^\r_1$.
\nl
This is clear.

Using (c), we see that 1.6(f) holds. Now the proof of 1.6(g) (in the geometric version) is entirely
similar to the proof of the corresponding statement in the $\ZZ$-graded case, 
see \cite{\GRAII, 3.14(c), 3.17}. Also, 1.6(h) follows immediately from 1.6(g).

We now assume that $m=1$. In this case 1.9(b) follows from 4.3(a) and 1.10(a) follows from 4.5(a).
Using (b) we see that 1.11(a) follows from the following statement: if $(\co,\ce)\in\ci_0$ and
$\dim\co=\sha(R)$, then $\co$ is the regular nilpotent orbit and $\ce=\bbq$; this is immediate.
The statement 1.12(a) follows from 6.8(b).

\subhead 8.4\endsubhead
In the remainder of this section when $m=\iy$ we write $\BB^\o,Z,\un{\un\YY}^\bul$ instead of ${}^1\BB^\o,{}^1Z_R,{}^1\un{\un\YY}^\bul$
(see B.6); note that $\BB^\o,Z,\un{\un\YY}^\bul$ are also defined when $m<\iy$.

We prove Theorem 0.3. 
For $\et\in\BB$ we denote by $\et^!$ the image of $b$ under the canonical bijection 
$\BB\lra Z$ (see 1.6 for $m<\iy$, B.6 for $m=\iy$). For $\et,\et'$ in $\BB$ we define $\MM'_{\et,\et'}\in\ZZ[v]$ by
$$\et'=\sum_{\et\in\BB}\MM'_{\et,\et'}\et^!.$$
Let $\BB\lra\ci_0$ be the bijection such that $\et\in\BB$ is mapped to $(\co,\ce)$
if $\et$ corresponds to $\tt_{\co,\ce}$ under the identification
$\VV=\QQ(v)\ot_\ca\ck_0$. By 3.6(b), this bijection has the property
stated in the theorem.

\subhead 8.5\endsubhead 
We prove Theorem 0.5. Let $\Th=\un{\un\YY}^\bul$. We define $\c:\BB@>>>\Th$ by 
$\et\m\o$ where $\et\in\BB^\o$ (see 1.6 when $m<\iy$ and B.6 when $m=\iy$).

We define $h':\Th@>\si>>G_0\bsl\fg_1^{nil}$
as in 8.3(a) if $m<\iy$ and as in 8.2(b) if $m=\iy$. With these definitions, the
theorem holds.

\subhead 8.6\endsubhead 
We prove Theorem 0.6. The set $\hW^{[\o]}$ and the bijection $\c\i(\o)\lra\hW^{[\o]}$ are defined in 1.13.
This proves the theorem.

\subhead 8.7\endsubhead 
The proof of Theorem 0.7 is contained in 1.12 (see 1.12(b)).

In the case where $\o\in\Th$ corresponds to a nilpotent orbit in $\fg$ which is not even, the set
$\hW^{[\o]}$ in 0.6 (with $m=1$) consists of certain irreducible representations of a proper subgroup
of $W$, hence it is not the same as $\hW^\o$ in 0.7 which consists of 
certain irreducible representations of $W$ itself.

\head Appendix A. An alternative definition of the PBW basis\endhead
\subhead A.1\endsubhead
In this appendix we assume that we are in the setup of 1.2; in particular we have $m<\iy$.
Let
$$\ovsc\YY=\YY-\cup_{\a\in R_1}\{y\in\YY;(y,\a)=1\},$$
$$\YY''=\YY-\cup_{N\in\ZZ-\{0\},\a\in R_{\bN}}\{y\in\YY;(y,\a)=N\}.$$
We have $\YY'\sub\YY''\sub\ovsc\YY$.

For $y,y'\in\ovsc\YY$ we say that $y\si y'$ if for any $\a\in R_1$ we 
have

$((y,\a)-1)((y',\a)-1)>0$.
\nl
This is an equivalence relation on $\ovsc\YY$. Let $\un{\ovsc\YY}$ 
be the set of equivalence classes (a finite set).
The following holds.

(a) {\it Let $c_1\in\un{\ovsc\YY}$, $c_2\in\un{\ovsc\YY}$ and let $y_1\in c_1\cap\YY'$,
$y'_1\in c_1\cap\YY'$, $y_2\in c_2\cap\YY'$, $y'_2\in c_2\cap\YY'$. Then $(y_1:y_2)=(y'_1:y'_2)$.}
\nl
This can be deduced from the arguments in the proof of
\cite{\LYII, 10.7(a)}; it can be also proved directly from the definitions.
It follows that for any $c\in\un{\ovsc\YY}$ there is a well defined element $T_c\in\VV$ such that
$T_c=I_\g$ for any $\g\in\un\YY'$ such that $\g\sub c$.

\mpb

We note that the definition of $\VV'$ in \cite{\LYII,\S10} is different from the one in this paper (it
is defined as a vector space with basis indexed by $\un{\ovsc\YY}$); however the vector space $\VV$
in \cite{\LYII,\S10} is the same as the one we use in this paper. The inner product on $\VV$ used
in \cite{\LYII} is of the form $x,x'\m s(x:\b(x'))$ where $(:)$ is as in this paper and $s\in\QQ(v)$
satisfies $s\in 1+v\ZZ[v]]$.

\subhead A.2\endsubhead
Note that if $w\in W_0$ and $c\in\un{\ovsc\YY}$ then $w(c)\in\un{\ovsc\YY}$. Thus $W_0$ acts naturally
on $\un{\ovsc\YY}$. For $c\in\un{\ovsc\YY}$ we denote by $R_{0,c}$ the set of roots in $R_0$
such that the corresponding reflection keeps $c$ stable. Let $\chR_{0,c}$ be the image of $R_{0,c}$
under $\chR\lra R$. Then $(\YY,\XX,(,),\chR_{0,c}\lra R_{0,c})$ is a root system. Let $W_{0,c}$ be 
the Weyl group of this root system viewed as a subgroup of $W_0$. 
Note that any $w\in W_{0,c}$ keeps $c$ stable.
Let $e_{W_{0,c}}=\sum_{w\in W_{0,c}}v^{2|w|}$ where $w@>>>|w|$ is the length function on $W_{0,c}$ for a
Coxeter group structure on $W_{0,c}$ determined by any choice of simple roots for $R_{0,c}$.
Let $\VV_\ca$ be the $\ca$-submodule of $\VV$ generated by $\{e_{W_{0,c}}\i T_c;c\in\un{\ovsc\YY}\}.$
We have the following result:

(a) {\it $\VV_\ca$ is equal to the $\ca$-submodule of $\VV$ generated by the canonical basis.}
\nl
Let $c\in\un{\ovsc\YY}$. For any $y\in\YY''\cap c$ let
$W_{0,y}$ be the subgroup of $W_0$ generated by reflections in the roots in
$R_0$ which are zero on $y$. We define $e_y\in\ZZ[v^2]$ in the same way as
$e_{W_{0,c}}$ above but replacing $W_{0,c}$ by $W_{0,y}$.
Let
$$\YY_0=\{x\in\YY;w(x)=x\qua \frl w\in W_{0,c}\}.$$
From the definition of $\ovsc\YY$ we see that $c$ has the following
``convexity'' property: if $y_1,y_2,\do,y_k$ are elements of $c$ then
$(y_1+y_2+\do+y_k)/k\in c$. Hence if $y'\in c$ then
$\sum_{w\in W_{0,c}}w(y')/\sha(W_{0,c})\in c$. We see that $\YY_0\cap c\ne\emp$.
Note that $c$ is open in $\YY$ and $\YY_0\cap c$ is open in $\YY_0$.
Also the affine hyperplanes in $\YY-\YY''$
do not contain $0$ hence their intersection with $\YY_0$ (which does contain
$0$) is a union of a discrete set of affine hyperplanes in $\YY_0$.
Hence $\YY''\cap\YY_0$ is dense in $\YY_0$. Since $\YY_0\cap c$ is open and
nonempty in $\YY_0$ it follows that $\YY''\cap\YY_0\cap c\ne\emp$.
Let $y_c\in\YY''\cap\YY_0\cap c$.
Then $y_c$ is fixed by the reflection with respect to any root in $R_{0,c}$
hence any such root is zero at $y_c$ so that $W_{0,y_c}=W_{0,c}$ and

(b) $e_{W_0,c}=e_{y_c}$. 
\nl
From \cite{\LYII, \S11} we see that the $\ca$-submodule of $\VV$
generated by
$$\{e_y\i T_c;c\in\un{\ovsc\YY},y\in \YY''\cap c\}$$
is equal to
the $\ca$-submodule of $\VV$ generated by the canonical basis of $\VV$.
Hence to prove (a) it is enough to show that for any $c\in\un{\ovsc\YY}$,
the $\ca$-submodule $M_c$ of $\QQ(v)$ generated by
$\{e_y\i;y\in\YY''\cap c\}$ is equal to $e_{W_{0,c}}\i\ca$.
For $y\in\YY''\cap c$ we have $W_{0,y}\sub W_{0,c}$. It follows that
$e_{W_{0,c}}/e_y\in \ZZ[v^2]$ so that $M_c\sub e_{W_0,c}\i\ca$.
From (b) we see that $e_{W_0,c}\i\ca\sub M_c$ hence $M_c=e_{W_0,c}\i\ca$.
This completes the proof of (a).

\subhead A.3\endsubhead
Let $\ti\BB'$ be the set of all $\et\in\VV_\ca$ such that $(\et:\et)\in1+v\ZZ[v]$. As in \cite{\LYII, 11.11} we see that
the following holds:

(a) {\it $\ti\BB'$ is a signed basis of the $\ca$-module $\VV_\ca$ (that is the union of a basis
with $(-1)$ times that basis; it is also a signed basis of the $\QQ(v)$-vector space $\VV$. There is a
unique $\ca$-basis $\ti\BB$ of $\VV_\ca$ such that for any $c\in\un{\ovsc\YY}$, the element $e_{W_0,c}\i T_c\in\VV_\ca$
is a $\NN$-linear combination of elements in $\ti\BB$.}
\nl
Under the identification $\VV=\QQ(v)\ot_\ca\ck_0$ (in 8.1), $\VV_\ca$ corresponds to $\ck_0$ and $\ti\BB$
corresponds to $\{\tt_{\co,\ce};(\co,\ce)\in\ci_0\}$. It follows that

(b) {\it $\ti\BB$ is the same as $\BB$ in 1.6(e).}

\subhead A.4\endsubhead
For any $m$-facet $\r$ we set
$$d_\r=\sha(\a\in R_0;(y,\a)<0)+\sha(\a\in R_1;(y,\a)\ge1)$$
where $y\in\r$. Let $[\r]$ be the set of all $\g'\in\un\YY'$ such that $\r$ is contained in the closure
of $\g'$. For any $\g'\in\un\YY'$ we write 
$I_{\g'}=\sum_{\et\in\ti\BB}N_{\et,\g'}\et\in\VV$  where $N_{\et,\g'}\in\ZZ[v,v\i]$.

 For $\et\in\ti\BB$ let $D(\et)\in\NN$ be the minimum of all integers $d_\r$ where $\r$ runs
through the $m$-facets such that $N_{\et,\g'}\ne0$ for some $\g'\in[\r]$.
For any $n\in\NN$ let $\ti\BB_n=\{\et\in\ti\BB;D(\et)=n\}$, $\ti\BB_{\le n}=\{\et\in\ti\BB;D(\et)\le n\}$.
Let $\VV_n$ (resp. $\VV_{\le n}$) be the $\QQ(v)$-subspace of $\VV$ spanned by $\ti\BB_n$ (resp. $\ti\BB_{\le n}$).
When $n=-1$ we set $\VV_{\le n}=0$.
We have the following result.

(a) {\it Let $n\in\NN$. There is a unique subspace $\VV^!_n$ of $\VV_{\le n}$ such that
$\VV_{\le n}=\VV_{\le n-1}\op\VV^!_n$ and $(\VV_{\le n-1}:\VV^!_n)=0$. Hence for any $\et\in\ti\BB_n$ there is a
unique element $\et^!\in\VV^!_n$ such that $\et-\et^!\in\VV_{\le n-1}$.}

\subhead A.5\endsubhead
The proof of A.4(a) is based on geometry. Let $\et\in\ti\BB$ and let $\tt_{\co,\ce}$
be the corresponding basis element of $\ck_0$.
Let $\r$ be an $m$-facet such that $N_{\et,\g'}\ne0$ for some $\g'\in[\r]$. 
Then $\g'=\ti\g=f_\r(\g)$ (notation of 1.4, 7.1) for some $\iy$-facet $\g$ of $\YY$ with respect to $R(\r)$.
Now some shift of $\ce^\sha$ is a direct summand of the complex $\tK^{\ti\g}$ which (as in 7.1) is  obtained 
from $\tK^\g$ on $\fl^\r_1$ by the induction procedure in \cite{\LYI, \S4} or 7.1.
The support of $\tK^{\ti\g}$ has dimension $\le \dim E''$ where $E''$ is as in 3.3.
 Hence $\dim\co\le\dim E''=\dim\fp^\r_1+\dim G_0-\dim P^\r=d_\r$ 
(notation of 3.3).
Now let $y\in C_\co$, $\r=\r_y$ be as in 5.4. Then by the results in 5.4,
we have $N_{\et,\g'}\ne0$ for some $\g'\in[\r]$. Moreover, by 8.3(b) we have $d_\r=\dim\co$. 
We see that $\dim\co=D(\et)$. Hence under the identification
$\VV=\QQ(v)\ot_\ca\ck_0$ (in 8.1), for $n\in\NN$, $\VV_n$ becomes the subspace
of $\QQ(v)\ot_\ca\ck_0$ spanned by $\{\tt_{\co',\ce'};(\co',\ce')\in\ci_0,\dim\co'=n\}$
and $\VV_{\le n}$ becomes the subspace of $\QQ(v)\ot_\ca\ck_0$ spanned by
$\{\tt_{\co',\ce'};(\co',\ce')\in\ci_0,\dim\co'\le n\}$.
In terms of the basis $\BB$ and its partition $\sqc_{\o\in\un{\un\YY}^\bul}\BB^\o$ in 1.6(e) we have
that $\VV_n$ is the subspace of $\VV$ spanned by $\sqc_{\o;d(\o)=n}\BB^\o$ and
$\VV_{\le n}$ is the subspace of $\VV$ spanned by $\sqc_{\o;d(\o)\le n}\BB^\o$. (We use 8.3(c).)
From 1.6(f)(g) we see that $\VV^!_n$ is the subspace of $\VV_{\le n}$ spanned by 
$\sqc_{\o;d(\o)=n}Z^\o_R$ and that if $\et\in\BB^\o$ with $d(\o)=n$ then $\et^!$ is the unique element
of $Z^\o_R$ such that $\et-\et^!\in\VV_{\le n-1}$. This completes the proof of A.4(a). 

\subhead A.6\endsubhead
Let $\ti Z$ be the subset of $\VV$ consisting of the elements $\et^!$ for various $n\in\NN$ and various
$\et\in\ti\BB_n$. Then $\ti Z$ is an $\ca$-basis of $\VV_\ca$ and an $\QQ(v)$-basis of $\VV$. It is in fact
equal to the basis $Z_R$ in 1.6. The present definition does not rely on the results in \S1, although the
proof of its correctness does.

\head Appendix B. $\ZZ$-graded root systems (by G. Lusztig)\endhead
\subhead B.1\endsubhead
In this appendix we reformulate the results in \cite{\GRAII} in a form which can be used in this paper.

Let $(\YY,\XX,(,),\chR\lra R)$ be as in 1.1. Let $\cs$ be the collection of linear hyperplanes
$$\{\{y\in\YY;(y,\a)=0\};\a\in R\}.$$
Now $\cs$ determines a set of facets called {\it $\iy$-facets} as follows.
For $y_1,y_2$ in $\YY$ we write $y_1\si y_2$ if for any $\a\in R$ we have 
$(y_1,\a)\ge0\Lra(y_2,\a)\ge0$. If $y_1\si y_2$, then for any $\a\in R$ we have 
$(y_1,\a)>0\Lra(y_2,\a)>0$. (Indeed, assume that $(y_1,\a)>0\Lra(y_2,\a)\not>0$. 
We have $(y_2,\a)=0=(y_2,-\a)=0$ hence 
$(y_1,-\a)\ge0$ and $(y_1,\a)\le0$, contradicting $(y_1,\a)>0$.) We deduce that if 
$y_1\si y_2$ then for any $\a\in R$ we have $(y_1,\a)=0\Lra(y_2,\a)=0$. Now $\si$ is 
an equivalence relation; the equivalence classes are the $\iy$-facets.
For example,
$$[0]:=\{y'\in\YY;(y',\a)=0\qua\frl\a\in R\}$$
is an $\iy$-facet of $R$ and 
$$\YY'=\YY-\cup_{\a\in R}\{y\in\YY;(y,\a)=0\}$$
is a union of $\iy$-facets called {\it $\iy$-alcoves}.
Let $\un\YY$ be the set of $\iy$-facets. Let $\un{\YY'}$ be the set of $\iy$-alcoves.
The $W$-action on $\YY$ induces a $W$-action on $\un\YY$ preserving $\un{\YY'}$.

For any $\r\in\un\YY$ let $R(\r)=\{\a\in R;(y,\a)=0\}$ where $y\in\r$;
this is independent of the choice of $y$. Let $\chR(\r)$
be the image of $R(\r)$ under $\chR\lra R$. Then $(\YY,\XX,(,),\chR(\r)\lra R(\r))$ is a 
root system. 

\subhead B.2\endsubhead
A {\it $\ZZ$-grading} of $R$ is a collection $R_*=(R_N)_{N\in\ZZ}$ where $R_N$ are subsets
of $R$ such that 
$R=\sqc_{N\in\ZZ}R_N$ and such that for some $y\in\YY$ we have 
$R_N=\{\a\in R;(y,\a)=N\}$ for all $N\in\ZZ$; we can assume that $y\in\la\chR\ra$; 
then $y$ is uniquely determined and is denoted by $y_{R_*}$.

\subhead B.3\endsubhead
We now fix a $\ZZ$-grading $R_*$ of $R$. We also fix $\d\in\{1,-1\}$.

Let $W_0$ be the subgroup of $W$ generated by the reflections with respect to roots in 
$R_0$. The obvious $W_0$-action on $R$ leaves stable each of the subsets $R_N,N\in\ZZ$.
Let $e_{W_0}=\sum_{w\in W_0}v^{2|w|}$ where $w@>>>|w|$ is the length function on $W_0$ 
for a Coxeter group structure on $W_0$ determined by any choice of simple roots for $R_0$. 
We have $\bb(e_{W_0})=v^{-\sha(R_0)}e_{W_0}$.

For $y,y'$ in $\YY'$ we define
$$\t(y,y')=\t(y',y)=\sha\{\a\in R_\d;(y,\a)(y',\a)<0\}-\sha\{\a\in R_0;(y,\a)(y',\a)<0\}\in\ZZ.$$
This is independent of $\d$ (we use that $\a\m-\a$ is a bijection $R_\d\lra R_{-\d}$); 
this justifies our notation. For $y,y'$ in $\YY'$ we define
$$(y:y')=e_{W_0}\sum_{w\in W_0}v^{\t(w(y),y')}\in\ca.$$
Let $\VV'=\VV'_R$ be the $\QQ(v)$-vector space with basis $\{I_\g;\g\in\un{\YY}'\}$. We define 
a bilinear form $(:):\VV'\T\VV'@>>>\QQ(v)$ by $(I_\g:I_\g')=(y:y')$ where $y\in\g$, 
$y'\in\g'$; this is independent of the choice of $y,y'$. 
This form is symmetric since $\t(y,w(y'))=\t(y',w\i(y))$ for $w\in W_0$. Let 
$\car=\{x\in\VV';(x:\VV')=0\}$, $\VV=\VV_R=\VV'/\car$. Then $(:)$ induces a symmetric nondegenerate
bilinear form on $\VV$ denoted again by $(:)$. 
For $\g\in\un{\YY'}$, the image in $\VV$ of $I_\g\in\VV'$ is denoted again by $I_\g$.

Define a semilinear involution 
$\b:\VV'@>>>\VV'$ by $\b(I_\g)=I_\g$ for all $\g\in\un{\YY}'$. Define a $\QQ(v)$-linear 
involution $\s:\VV'@>>>\VV'$ by $\s(I_\g)=I_{-\g}$ for all $\g\in\un{\YY'}$. We show:

(a) {\it For $\x,\x'$ in $\VV'$ we have $\bb((\b(\x):\b(\x')))=v^{-\sha(R_\d)}(\x:\s(\x'))$.}
\nl
We can assume that $\x=I_\g,\x'=I_{\g'}$ with $\g,\g'$ in $\un{\YY'}$. We must show:
$$v^{-\sha(R_0)}e_{W_0}\sum_{w\in W_0}v^{-\t(\g,w(\g'))}
=v^{-\sha(R_\d)}e_{W_0}\sum_{w\in W_0}v^{\t(\g,-w(\g'))}.$$
Let $w\in W_0$ and let $y\in\g,y'\in w(\g')$. It
is enough to show:
$$\align&\sha(R_0)-\sha(R_\d)+\sha(\a\in R_\d;(y,\a)(y',\a)<0)-\sha(\a\in R_0;(y,\a)(y',\a)
<0)\\&+\sha(\a\in R_\d;(y,\a)(y',\a)>0)-\sha(\a\in R_0;(y,\a)(y',\a)>0)=0;\endalign$$
this is clear. This proves (a).

\mpb

From (a) we see that $\b(\car)\sub\car$ hence $\b$ induces a semilinear involution 
$\VV@>>>\VV$ denoted again by $\b$.

\subhead B.4\endsubhead
Let $\r\in\un\YY$. Then $R(\r)$ has a $\ZZ$-grading $R(\r)_*$
where $R(\r)_N=R(\r)\cap R_N$ for all $N$. Hence $y_{R(\r)_*}\in\la\chR(\r)\ra\sub\YY$ is 
defined. We denote by $\YY'_\r$, $\un{\YY}'_\r$, $[0]_\r$ the analogues of $\YY'$, $\un{\YY}'$,
$[0]$ when $R$ is replaced by $R(\r)$.

We define a map 
$f_\r:\un{\YY}'_\r@>>>\un{\YY}'$ as follows. Let $\g\in\un{\YY}'_\r$ and let $y_1\in\g$, 
$y\in\r$. We have $(y_1,\a)\ne0$ for any $\a\in R(\r)$. We can assume that
$(y,\a)\in\ZZ,-1<(y_1,\a)<1$ for any $\a\in R$. We show:

(a) $y+y_1\in\YY'$.
\nl
If $\a\in R(\r)$ then $(y,\a)=0$ and $(y_1,\a)\ne0$ hence $(y+y_1,\a)\ne0$. If $\a\n R(\r)$
then $(y,\a)\in\ZZ-\{0\}$ and $-1<(y_1,\a)<1$ so that $(y+y_1,\a)\ne0$. This proves (a).

Now let $y'_1\in\g$, $y'\in\r$ be such that $(y',\a)\in\ZZ,-1<(y_1,\a)<1$ for any $\a\in R$. 
By (a) we have $y'+y'_1\in\YY'$. We show:

(b) $y+y_1\si y'+y'_1$.
\nl
Assume that for some $\a\in R$, $(y+y_1,\a),(y'+y'_1,\a)$ have different 
signs. If $\a\in R(\r)$ then $(y,\a)=(y',\a)=0$, so that $(y_1,\a)$, $(y'_1,\a)$
have different signs; 
this contradicts $y_1\in\g,y'_1\in\g$. If $\a\n R(\r)$ then $|(y,\a)|\ge1$; since $-1<(y_1,\a)<1$, we 
see that $(y,\a)$ has the same sign as $(y+y_1,\a)$. Similarly $(y',\a)$ has the same sign as 
$(y'+y'_1,\a)$. Thus $(y,\a),(y',\a)$ have different signs. This contradicts $y\si y'$ and proves (b).

We see that $\g\m y+y_1$ is a well-defined map
$f_\r:\un{\YY}'_\r@>>>\un{\YY}'$. 

\subhead B.5\endsubhead
Let $\VV'_{R(\r)},\VV_{R(\r)},(:)_\r$ be the analogues of $\VV',\VV,(:)$ where $R_*$ is replaced by 
$R(\r)_*$. We define a linear map \footnote{This corresponds to $\ZZ$-graded parabolic induction.}
$\VV'_{R(\r)}@>>>\VV'$ by sending the basis element indexed 
by $\g\in\un{\YY}'_\r$ to $I_{f_\r(\g)}$. By an argument using geometry in \cite{\GRAII} (see also 8.2)
one can show

(a) {\it this maps the radical of $(:)_\r$ on $\VV'_{R(\r)}$ into the radical of $(:)$ on $\VV'$ 
hence it induces a linear map $\VV_{R(\r)}@>>>\VV$ denoted again by $f_\r$.}

\subhead B.6\endsubhead
Using induction on $\sha(R)$ we define a subset ${}^\d\un\YY^\bul$ of $\un\YY$ and for each 
$\r\in{}^\d\un\YY^\bul$ we define a nonempty subset ${}^\d Z^\r_R$ of $\VV$. 
Assume first that $R=\emp$. Then $\VV'=\VV$ is one dimensional with basis $\{I_{[0]}\}$. We
define  ${}^\d\un\YY^\bul=\{[0]\}$, ${}^d Z^{[0]}_R=\{I_{[0]}\}$. 
Next we assume that $R\ne\emp$. 
Let $\r\in\un\YY$. Assume first that $\r\ne[0]$. We have $\sha(R(\r))<\sha(R)$. 
We say that $\r\in{}^\d\un\YY^\bul$ if 
$(y_{R(\r)_*}-y_{R_*})\d\in\r$ and $[0]_\r\in{}^\d\un\YY_\r^\bul$
(which is already defined); we set ${}^\d Z^\r_R=f_\r({}^\d Z^{[0]_\r}_{R(\r)})$.
It remains to decide whether $[0]$ is in ${}^\d\un\YY^\bul$ or not and, if it is, to define
${}^\d Z^{[0]}_R$.
Let ${}^\d Z'=\cup_{\r\in{}^\d\un{\un\YY}^\bul;\r\ne[0]}{}^\d Z^\r_R\sub\VV$. 
Let $\cl$ be the $\ZZ[v]$-submodule of $\VV$ generated by ${}^1Z'\cup{}^{-1}Z'$, let
$\p:\cl@>>>\cl/v\cl$ be the obvious map and let ${}^1\z'=\p({}^1Z'),{}^{-1}\z'=\p({}^{-1}Z')$,
$\z={}^1\z'\cup{}^{-1}\z'$.
By an argument using geometry in \cite{\GRAII} (see also 8.2) one can show: 

(a) {\it $\z$ is a $\ZZ$-basis of $\cl/v\cl$. For any $b\in\z$ there is a unique 
element $\tb\in\cl$ such that $\p(\tb)=b$ and $\b(\tb)=\tb$. Moreover,
$\{\tb;b\in\z\}$ is a $\ZZ[v]$-basis $\BB$ of $\cl$ and a $\QQ(v)$-basis of $\VV$.}

If $\z={}^\d\z'$ then we declare $[0]\n{}^\d\un\YY^\bul$. 
If $\z\ne{}^\d\z'$ then we declare $[0]\in{}^\d\un\YY^\bul$.
By an argument using geometry in \cite{\GRAII} (see also 8.2) one can show: 

(b) {\it Assume that $\z\ne{}^\d\z'$; let ${}^\d\LL'$ be the $\QQ(v)$-vector subspace of $\VV$ 
spanned by ${}^\d Z'$. Then for any $b\in\z-{}^\d\z'$ there is a unique element $\hat b\in\VV$ such that
$\tb-\hat b\in{}^\d\LL'$, $(\hat b:{}^\d\LL')=0$. The map $b\m\hat b$ is
a bijection of $\z-{}^\d\z'$ onto a subset of $\cl$ denoted by ${}^\d Z^{[0]}$.}
\nl
This completes the inductive definition of ${}^\d\un\YY^\bul$ and of the subsets
${}^\d Z^\r_R$ (for $\r\in{}^\d\un\YY^\bul$).

Note that if $\r=[0]$ is in ${}^\d\un\YY^\bul$ then it is again true that $(y_{R(\r)_*}-y_{R_*})\d\in\r$. 
Indeed, in this case we have $R(\r)=R$ so that $(y_{R(\r)_*}-y_{R_*})\d=0$.

Now the $W_0$-action on $\un\YY$ restricts to a $W_0$-action on 
${}^\d\un\YY^\bul$. Let ${}^\d\un{\un\YY}^\bul$ be the set of orbits of this last action.
\footnote{This is a combinatorial version of the set of $G_0$-orbits on $\fg_1$, see 8.2.}
Note that if $\r,\r'$ in ${}^\d\un\YY^\bul$ are in the same $W_0$-orbit then
${}^\d Z^\r_R={}^\d Z^{\r'}_R$. Hence for any $\o\in{}^\d\un{\un\YY}^\bul$ we can define
${}^\d Z^\o_R={}^\d Z^\r_R$ where $\r$ is any $\iy$-facet in $\o$.
Let ${}^\d Z_R=\cup_{\o\in{}^\d\un{\un\YY}^\bul}{}^\d Z^\o_R$.
By an argument using geometry in \cite{\GRAII} one can show: 

(c) {\it We have ${}^\d Z_R=\sqc_{\o\in{}^\d\un{\un\YY}^\bul}{}^\d Z^\o_R$.
If $\o\in{}^\d\un{\un\YY}^\bul$, $\r\in\o$, then $f_\r:{}^\d Z^{[0]_\r}_{R(\r)}@>>>{}^\d Z^\o_R$ is a
bijection independent of the choice of $\r$. Moreover, ${}^\d Z_R$ is a 
$\ZZ[v]$-basis of $\cl$ and a $\QQ(v)$-basis of $\VV$ which we call a PBW basis. The map 
$\x\m\p(\x)$ defines a bijection ${}^\d Z_R@>\si>>\z$. Hence there is a unique bijection 
${}^\d Z_R@>\si>>\BB$, $\x\m\un\x$ defined by the requirement that $\p(\x)=\p(\un\x)$ for any 
$\x\in{}^\d Z_R$. For $\o\in{}^\d\un{\un\YY}^\bul$ let ${}^\d\BB^\o$ be the subset of $\BB$
corresponding to ${}^\d Z_R$ under this bijection. We have
$\BB=\sqc_{\o\in{}^\d\un{\un\YY}^\bul}{}^\d\BB^\o$.}
\nl
Thus $\VV$ has two PBW bases: ${}^1Z_R$ and ${}^{-1}Z_R$ and one basis $\BB$ which we call
{\it canonical basis}. 

\subhead B.7\endsubhead
We say that $R_*$ in B.3 is {\it rigid} if ${}^1Z^{[0]}_R\ne\emp$ or equivalently if
${}^{-1}Z^{[0]}_R\ne\emp$. (The equivalence follows by an argument using geometry in
\cite{\GRAII}.)

Let $\r\in\un\YY$. We say that $\r$ is {\it $\d$-rigid} if $\r\in{}^\d\un\YY^\bul$.

\mpb

{\it Errata to \cite{\LYI}.}

page 277, line 3 of Contents. Replace $\ZZ/\m$ by $\ZZ/m$.

page 280, line 5. Replace ``...for large $m$, a $\ZZ/m$-grading is the same as a $\ZZ$-grading, so 
that in this case ...results of \cite{L4}.''
by:  ``...a $\ZZ$-grading can be viewed as a $\ZZ/m$-grading for large $m$. (The converse does
not hold.)''

page 287. Title of Section 2. Replace $\ZZ/\m$ by $\ZZ/m$.

page 303, line 3 of 5.1. Relace $\fp''$ by $\fp''_*$.

page 304, line 7. Replace $\dot\fg_\d=...''$ by 
$\dot\fg_\d=\{(gP_0,z)\in G_0/P_0\T \fg_\fd;\Ad(g\i)z\in\p\i(\fl_\et^0)\}$.

page 305, lines 11-15. Replace ${}_{k,k'}\fg_{\bN}$ by ${}_{k',k''}\fg_{\bN}$.

\enddocument